\documentclass{amsart}
\usepackage{graphicx}
\usepackage{graphicx,esint}
\usepackage{amssymb,amscd,amsthm,amsxtra}
\usepackage{latexsym}
\usepackage{epsfig}
\newtheorem{thm}{Theorem}[section]
\newtheorem{cor}[thm]{Corollary}
\newtheorem{lem}[thm]{Lemma}
\newtheorem{prop}[thm]{Proposition}
\theoremstyle{definition}
\newtheorem{defn}[thm]{Definition}
\theoremstyle{remark}

\numberwithin{equation}{section}
\newcommand{\R}{\mathbb R}

\newcommand{\eps}{\epsilon}
\newcommand{\ov}{\overline}
\newcommand{\p}{\partial}
\renewcommand{\theenumi}{\roman{enumi}}
\newcommand{\comment}[1]{}
\begin{document}

\title[Monotone free boundaries]{Existence and regularity of monotone solutions to a free boundary problem}
\author{Daniela De Silva}
%\author{}%
%\address{}%
%\email{}%

%\thanks{}%
%\subjclass{}%
%\keywords{}%

%\date{}%
%\dedicatory{}%
%\commby{}%
% ----------------------------------------------------------------
\begin{abstract} We investigate existence and regularity properties
of one-phase free boundary graphs, in connection with the question
of  whether there exists a complete non-planar free boundary graph
in high dimensions.

\end{abstract}
\maketitle

\section{Introduction}

Let $\Omega$ be an open connected subset of $\R^n$, and consider
the energy functional
\[J(u,\Omega)= \int_{\Omega}
(|\nabla u|^2 + \chi_{\{u>0\}}).\] In \cite{AC}, Alt and
Caffarelli analyzed the question of the existence and regularity
of a minimizer $u$ of $J(\cdot,\Omega)$. They developed a partial
regularity theory for the free boundary of $u$, that is
\[F(u)=\p \{u>0\}
\cap \Omega,\] showing that its reduced part $F^*(u)$ is locally
$C^{1,\alpha}.$ Higher regularity results of Kinderhlerer and
Nirenberg \cite{KN}, then imply that $F^*(u)$ is locally analytic.
In \cite{AC} the authors also proved that in two dimensions,
$F(u)$ does not have singularities. Subsequently in \cite{W2},
Weiss showed that there exists a critical dimension $k $, $3 \leq
k \leq +\infty$, such that energy minimizing free boundaries are
smooth for $n < k$. These results draw on a strong analogy with
the theory of minimal surfaces, for which it is known that the
critical dimension is 8. In \cite{CJK}, Caffarelli, Jerison and
Kenig proved that there are no singular free boundary minimizers
in dimension $n=3$, which yields $k \geq 4$. Their proof suggests
that $k=7$, but the problem remains still open. In \cite{DJ} De
Silva and Jerison showed that $k \leq 7$, by providing the first
example of a singular energy minimizing free boundary in dimension
$n=7$. Analogously, for the theory of minimal surfaces, the Simons
cone, provides an example of a singular set of minimal perimeter
in dimension $n=8$.

The purpose of this note is to pursue even further the analogy
between the theory of minimal surfaces and free boundary
regularity, precisely we turn to free boundary graphs. Our
motivation lies in the question of whether there exists a complete
non-planar free boundary graph, i.e a classical solution on $\R^n$
to the problem, \vspace{0.5mm}
\begin{equation}\label{exgraph}
  \begin{cases}
    \Delta u = 0 & \text{in $\{u>0\}$}, \\
    |\nabla u|=1 & \text{on $\partial \{u>0\}$},\\
    \partial \{u>0\}   & \text{is a non-planar graph in the $x_{n}$ direction.}
  \end{cases}
\end{equation}

\vspace{2mm}

\noindent This is the analogue of the celebrated Bernstein problem
for minimal graphs. Precisely, it is known that planes are the
only complete smooth minimal graphs in $\R^n$, when $n\leq 8$ (see
for example \cite{G}). This result is sharp since in \cite{BDG},
Bombieri, De Giorgi and Giusti proved the existence of a non
affine minimal graph in dimension 9, which turns out to be
strictly related to the existence of the Simons cone, one
dimension lower. The result in \cite{DJ}, then naturally raises
the analogous question for free boundary graphs. In analogy with
the minimal surfaces theory, we expect that a global smooth
solution to (\ref{exgraph}) exists in dimension $8$ or higher.

The first step towards constructing such an example, is to develop
a local theory which is the analogue of the existence and
regularity theory for the minimal surface equation in a ball, when
the boundary data is smooth. Then, a limiting argument provides a
global solution. In order to prevent this global solution from
being planar, one wishes to control its behavior by trapping it
between a given global subsolution and a given global
supersolution. However, since ordinary (strong) comparison results
are not available in the free boundary context (see Lemma
\ref{complemma}), that is if two solutions are one greater than
the other one on the boundary then the inequality may not be
preserved in the interior, this ``trapping" is not
straightforward. Thus, we develop our local theory by constructing
a solution which is trapped in between a given subsolution and
supersolution. Moreover, since we wish to preserve the graph
property in the limit, we also need to construct a local solution
which enjoys a certain regularity property (density property),
which for example, would not be guaranteed if we were to construct
our solution via a standard Perron method. For this reason, we
seek a solution which is also an energy minimizer to the energy
functional $J$ among a certain class of competitors.

\vspace{1.5mm}

Our result is the following. Let $R, h_R>0$, and let $C_R$  denote
the cylinder $C_R=\mathcal{B}_R(0) \times \{|x_{n}|<h_R\},$ and
$S_R$ denotes the sides of the cylinder $C_R$, $S_R=\partial
\mathcal{B}_R(0) \times \{|x_{n}|\leq h_R\}$; consider the
one-phase free boundary problem,

\begin{equation}
  \begin{cases}\label{FB1intro}
    \Delta u = 0 & \text{in $C_R^+(u):=\{x \in C_R : u(x)>0\}$}, \\
   |\nabla u|  = 1 & \text{on $F(u):=(\partial C^+_R(u)) \cap C_R$}.
  \end{cases}
\end{equation}

\

\begin{thm}\label{localintro} Assume that, there exist a strict
smooth subsolution $V_1$ and a strict smooth supersolution $V_2$
to $(\ref{exgraph})$ in $\R^{n}$, such that
\begin{enumerate}
\item $0 \leq V_1 \leq V_2$ on $\R^{n}; 0 \in \{V_2 >0\} \cap \{V_1=0\}^\circ;$
\item $\p_{n} V_i >0$ in $\ov{\{V_i>0\}},$ for $i=1,2$.
\end{enumerate}
Then, for each $R>0$ and $h_R$ sufficiently large, there exists
$u_R$ viscosity solution to $(\ref{FB1intro})$, such that $u_R$ is
minimizes $J$ over $K_R := \{v \in H^1(C_R) \ | \ V_1\leq v \leq
V_2 \ \ \text{a.e on $C_R$}, \ \ v = V_2 \ \ \text{on $S_R$}\},$
$u_R$ is monotone increasing in $\ov{C_R^+(u)}$ in the $x_{n}$
direction, and $V_1 \leq u_R \leq V_2$. Moreover, in the interior,
$F(u_R)$ is a smooth graph in the $x_{n}$ direction.
\end{thm}

For the precise definition of viscosity solution, we refer the
reader to Section 2. We remark that the proof of the result in
\cite{DJ}, provides a clear indication of how to construct
functions $V_1$ and $V_2$ satisfying the assumptions above, when
$n \geq 8.$ We plan on constructing these functions in a future
paper.

The main tools to achieve the existence part in our Theorem are
blow-up and domain variation techniques. Then, the fact that $u_R$
is also a viscosity solution, allows us to use maximum principle
techniques and a continuity argument, to compare $u_R$ with a
family of subsolutions, which are suprema of vertical translates
of $u_R$ over balls (supconvolutions). This yields the desired
Lipschitz behavior (hence smoothness) of the free boundary of
$u_R$.

As already observed, the second step towards constructing a global
solution to (\ref{exgraph}) is a limiting argument as $R
\rightarrow +\infty$. In the theory of minimal surfaces, the
convergence to a global solution is guaranteed by a very powerful
tool, that is the a-priori estimate of the gradient of a solution
to the minimal surface equation. In a forthcoming paper, we prove
the analogue of such a tool in the free boundary context
\cite{DJ2}.

Here, in order to preserve the graph property in the limit, we
prove that the positive phase of our solution is a
non-tangentially accessible (NTA) domain, that is, it enjoys a
certain scale-invariant connectivity property (see Section 2 for
the precise definition of NTA domains.) Then, a limiting argument
allows us to prove the following:

\begin{thm} \label{z} Assume that, there exist a strict smooth subsolution $V_1$ and a
strict smooth supersolution $V_2$ to $(\ref{exgraph})$ in $\R^n$,
such that:
\begin{enumerate}
\item $0 \leq V_1 \leq V_2$ on $\R^n$;
\item $\p_{n} V_i > 0$ on $\ov{\{V_i > 0\}}, i=1,2;$
\item $\lim_{r \rightarrow \infty}{\displaystyle\frac{V_1(rx)}{r}}\geq U(x').$
\end{enumerate}
Then, there exists a global energy minimizing viscosity solution
$u$ to:\[\Delta u=0 \ \ \ \text{in}\ \ \ \{u>0\},  \ \ \ |\nabla
u|=1 \ \ \ \text{on} \ \ \ F(u),\] such that $u$ is monotone
increasing in $\ov{\{u>0\}}$ in the $x_{n}$ direction, and $F(u)$
is a continuous non-planar graph, with a universal modulus of
continuity on each compact subset of $\R^n$. Moreover, $F(u)$ is
NTA.
\end{thm}

Here $U$ is the singular global minimizer in $\R^{n-1}$ from
\cite{DJ}, interpreted as a function of $n$ variables. Hypothesis
(iii) is used only to prevent $F(u)$ from being planar. While we
could weaken this assumption, its motivation lies in the fact
that, in analogy with the minimal surfaces theory, we expect a
smooth non-affine free boundary graph $u$ to blow down to an
energy minimizing solution.

The NTA property of $F(u)$ is proved by the means of a
monotonicity formula \cite{ACF} for $\nabla u$, together with
non-degeneracy properties of $u$. The proof follows from arguments
in \cite{ACS}. Then, exploiting the known behavior of positive
harmonic functions in NTA domains \cite{JK}, we derive that $F(u)$
cannot contain vertical segments.

The paper is organized as follows. In Section 2, we introduce some
notation and definitions. In Section 3, we show the existence of a
``trapped" monotone energy minimizing viscosity solution to our
one-phase free boundary problem on a cylinder. Then, in Section 4,
we show that such a solution is indeed smooth in the interior. In
Section 5, we prove that the free boundary of our solution is
locally NTA.  Finally, in Section 6, we prove the existence of a
global monotone viscosity solution whose free boundary is a graph
in the vertical direction, with a universal modulus of continuity
on each compact of $\R^n$, trapped between two given graphs.

\section{Notation and definitions.}

In this section we collect some notation and definitions which
will be used throughout this paper.

 A point $x \in \R^n $ will be
denoted by $(x', x_n)$, with $x'=(x_1,...,x_{n-1}).$ A ball of
radius $r$ in $\R^{n-1}$, will be denoted by $\mathcal{B}_r$,
while a ball of radius $r$ in $\R^n$, will be denoted by $B_r.$
When specifying the center $x$ of the ball, we will use either
$B_r(x)$ or $B(x,r).$ Also, $\Omega$ denotes an open bounded
connected subset of $\R^n.$

%Let $V_1 \leq V_2$ be non-negative functions on $\R^n$, satisfying
%(\ref{graph}), and let $0 \in \{V_2>0\} \cap \{V_1=0\}^\circ$. We
%will denote by \[C_R = \mathcal{B}_R(0) \times \{|x_n|<h_R\},\]
%with \[h_R \geq \max\{2d_R(V_1),2d_R(V_2), R\}.\] Roughly
%speaking, we consider a cylinder tall enough, so that the free
%boundaries of $V_1$ and $V_2$ ``exit" from the sides of the
%cylinder.

%Consider the energy functional,
%\[J(u)=\int (|\nabla u|^2 +
%\chi_{\{u>0\}}).\]

%\begin{defn}$u \in H^1_{loc}(\R^n)$ is a global
%minimizer for $J$, if and only if, for any ball $B \subset \R^n$,
%and any function $v \in H^1(B)$, such that $u-v \in H^1_0(B)$,
%$J(v,B) \geq J(u,B).$
%\end{defn}

%We introduce a notion of viscosity solution to the Euler equation
%for the functional $J$ (see \cite{C1}). First,

For any non-negative function $u$ on $\Omega$, set
$$\Omega^+(u)=\{x \in \Omega : u(x)>0\}, \ \ \Omega^- = \{x \in
\Omega :u =0 \}^\circ, \ \ \ F(u)= (\p \Omega^+(u))\cap \Omega.$$

Consider the one-phase free-boundary problem:
\begin{equation}
  \begin{cases}\label{FB1}
    \Delta u = 0 & \text{in $\Omega^+(u)$}, \\
   |\nabla u|  = 1 & \text{on $F(u)$}.
  \end{cases}
\end{equation}

We recall the following standard definition (see for example
\cite{C1}.)

%\begin{defn} \label{defsol}Let $u$ be a nonnegative continuous function in
%$\Omega$. We say that $u$ is a viscosity solution to (\ref{FB1})
%in $\Omega$, if and only if the following conditions are
%satisfied:
%\begin{enumerate}
%\item $\Delta u = 0$ in $\Omega^+(u)$;

%\item If $x_0 \in F(u)$ and $F(u)$ has at $x_0$ a one-sided tangent
%ball (i.e. there exists $B_{\epsilon}$ such that $x_0 \in
%\partial B_{\epsilon}$ and $B_{\epsilon}$ is contained either in
%$\Omega^+$ or in $\Omega^-$), then, for $\nu $ the unit radial
%direction of $\partial B_{\epsilon}$ at $x_0$ into $\Omega^+(u)$,
%\begin{equation}\nonumber u(x)=(x-x_0,\nu)^+ + o(|x-x_0|), \ \text{as $x\rightarrow x_0.$}
%\end{equation}
%\end{enumerate}
%\end{defn}

\begin{defn}\label{defsub} Let $v$ be a nonnegative continuous function in
$\Omega$. We say that $v$ is a viscosity subsolution (resp.
supersolution) to (\ref{FB1}) in $\Omega$, if and only if the
following conditions are satisfied:
\begin{enumerate}
\item $\Delta v \geq 0$ (resp. $\leq 0$) in $\Omega^+(v)$;

\item If $x_0 \in F(v)$ and $F(v)$ has at $x_0$ a
tangent ball $B_{\epsilon}$ from the positive (resp. zero) side
(i.e. $B_{\epsilon} \subset \Omega^+(v)$ (resp. $\Omega^-(v)$),
$x_0 \in
\partial B_{\epsilon}),$ then, for some $\alpha \geq 1$
(resp. $\alpha \leq 1$) and $\nu $ the unit inner (resp. outer)
radial direction of $\partial B_{\epsilon}$ at $x_0$,
\begin{equation}\nonumber
v(x)=\alpha (x-x_0,\nu)^+ + o(|x-x_0|), \ \text{as $x\rightarrow
x_0.$}
\end{equation}
\end{enumerate}
\end{defn}

\vspace{2mm}

If the constant $\alpha$ in Definition \ref{defsub} is strictly
greater (resp. smaller) than $1$, then $u$ is called a strict
viscosity subsolution (resp. supersolution).

When $u$ is simultaneously a viscosity subsolution and a viscosity
supersolution, then $u$ is called a viscosity solution.

\vspace{1.5mm}

Our viscosity solution will satisfy certain regularity properties,
which we now define. Let $u$ be a continuous non-negative function
on $\Omega$. Set,
\[d(x)=\text{dist}(x,F(u)).\]

\begin{defn}\label{nond} We say that $u$ is non-degenerate, if and only if,
for every $G \Subset \Omega$, there exists a constant $K=K(G)$
such that \[u(x) \geq K d(x),\] for all $x \in G^+(u),$ with
$B_{d(x)}(x) \subset G.$
\end{defn}

\begin{defn}\label{Inond} We say that  $u$ is $(I)$
non-degenerate, if and only if, for every $G \Subset \Omega$,
there exists a constant $K = K(G)$ such that, for any ball $B_r
\subset G$ centered at a free boundary point,\[\fint_{B_r}{u} \geq
K r.\]
\end{defn}

\begin{defn}\label{dens2} We say that $F(u)$ satisfies the density property $(D)$ if and only if:
\begin{enumerate}\renewcommand{\theenumi}{$D$}
\item for any $G \Subset \Omega,$ there exists a constant $c=c(G)<1,$
such that, for any ball $B_r \subset G$ centered at a free
boundary point, \[c \leq \frac{|B_r \cap \{u>0\}|}{|B_r|} \leq
1-c.\]
\end{enumerate}
\end{defn}

\vspace{2mm}

Such properties are crucial to use so-called ``blow-up"
techniques. Specifically, let $u$ be a non-negative, Lipschitz
continuous function in $\Omega$. Let $x_0 \in F(u)$, and let
$B_{r_k}(x_0) \subset \Omega$ be a sequence of balls with $r_k
\rightarrow 0$, as $k \rightarrow +\infty$. Consider the blow-up
sequence:\[u_k(x)= \frac{1}{r_k}u(x_0 + r_k x).\]

\vspace{1mm}

 \noindent Since for a given $D \Subset \R^n$ and
large $k$ the functions $u_k$ are uniformly Lipschitz continuous
in $D$, there exists a function $u_0: \R^n \longrightarrow \R$,
such that: \vspace{1mm}
\begin{align*}
& u_k \rightarrow u_0 \ \ \text{in} \ \ C^{0,\alpha}_{loc}(\R^n),
\ \ \text{for all $0 <
\alpha <1$};\\
& \nabla u_k \rightarrow \nabla u_0  \ \ \text{weakly star in
$L^{\infty}_{loc}(\R^n)$.}
\end{align*}
Moreover, $u_0$ is Lipschitz continuous in the entire space. The
globally defined function $u_0$ is called a blow-up of $u$. Using
the same argument as in \cite{F}, (see Chapter 3, Lemma 3.6), one
can prove the following.
\begin{lem} \label{blowup}Let $u$ be a non-negative function in $\Omega$, harmonic
in $\Omega^+(u)$, Lipschitz continuous and $(I)$ non-degenerate.
Assume that $u$ satisfies the density property $(D)$. Then the
following properties hold:
\begin{enumerate}\renewcommand{\theenumi}{\alph{enumi}}
\item $\p \{u_k >0 \} \rightarrow \p \{u_0 >0 \}$ in the Hausdorff distance;
\item $\chi_{\{u_k > 0\}} \rightarrow \chi_{\{u_0 > 0\}}$ in $L^1_{loc}$;
\item $\nabla u_k \rightarrow \nabla u_0$ a.e.
\end{enumerate}
Moreover, $u_0$ is $(I)$ non degenerate, and it satisfies the
density property $(D).$
\end{lem}

We also need the following result from \cite{W1}, which
characterizes blow-up limits of certain solutions (variational
solutions) to the free boundary problem \eqref{FB1}.

First we recall the following definition. Let,
\begin{equation}\label{J}J(u, \Omega)=\int_{\Omega} (|\nabla u|^2
+ \chi_{\{u>0\}}).\end{equation}

\begin{defn}
We define $u \in H^1_{loc}(\Omega)$ to be a variational solution
to (\ref{FB1}), if $u \in C(\Omega) \cap C^2(\Omega^+(u))$ and \[0
= -\frac{d}{d\eps}J(u(x+\eps \eta(x)))|_{\eps=0} =
\int_{\Omega}{(|\nabla u|^2 \mathrm{div} \eta - 2\nabla u D\eta
\nabla u + \chi_{\{u>0\}} \mathrm{div} \eta)}\] for any $\eta \in
C_0^1(\Omega,\R^n)$.
\end{defn}

\begin{lem}\label{weiss} Let $u$ be a variational solution to
$(\ref{FB1})$ in $\Omega$, and assume that $u$ is Lipschitz
continuous and satisfies the density property $(D)$. Then any blow
up limit of $u$ is homogeneous of degree 1.
\end{lem}

Such blow-up technique will allow us to prove that our energy
minimizing solution is indeed a trapped viscosity solution, which
also enjoys the desired regularity properties . These regularity
properties in Definitions
\eqref{nond}-\eqref{Inond}-\eqref{dens2}, are crucial in the proof
that the free boundary of our solution is locally NTA.

We conclude this section by recalling the notion of
Non-Tangentially Accessible (NTA) domains (see \cite{JK}). Let $D$
be a bounded domain in $\R^n$. A $M$-non-tangential ball in D, is
a ball $B_r \subset D$, such that: $Mr > \text{dist}(B_r,\p D)
> M^{-1}r.$ For $x_1,x_2 \in D,$ a Harnack chain from $x_1$ to $x_2$ in $D$ is
a sequence of $M-$non-tangential balls, such that the first ball
contains $x_1,$ the last contains $x_2$, and such that consecutive
balls intersect.

\begin{defn}A bounded domain $D$ in $\R^n$ is called NTA, when there exist
constants $M$ and $r_0 >0$ such that:
\begin{enumerate}
\item Corkscrew condition. For any $x \in \p D,$ $r < r_0,$ there exists $y=y_r(x) \in D$ such
that $M^{-1} r < |y-x|<r$ and $\text{dist}(y,\p D) > M^{-1}r;$
\item The Lebesgue density of $D^{c}$ at any of its points is bounded below uniformly by a positive constant $C$, i.e for
all $x \in \partial D, 0 < r <r_0,$ $$\frac{|B(x,r) \setminus D
|}{|B(x,r)|} \geq C;$$
\item Harnack chain condition. If $\eps >0$ and $x_1, x_2$ belong to $D$, $\text{dist}(x_j,\p
D)>\eps$ and $|x_1- x_2|< C\eps,$ then there exists a Harnack
chain from $x_1$ to $x_2$ whose length depends on $C$, but not on
$\eps$.
\end{enumerate}
\end{defn}

\section{Local existence theory.}

In this section we prove the existence of an energy minimizing
viscosity solution $u_R$ as in Theorem \ref{localintro}. We also
show that $u_R$ enjoys the regularity properties from Definitions
\eqref{nond}-\eqref{Inond}-\eqref{dens2}. The regularity of the
free boundary $F(u_R)$ will be investigated in the next section.

Precisely, for $R, h_R>0,$ let
$$C_R := \mathcal{B}_R(0) \times \{|x_n|<h_R\}, \ \ \ S_R= \p
\mathcal{B}_R(0)\times \{|x_n| \leq h_R\}.$$

We prove the following result.

\begin{thm}\label{localex} Assume that, there exist a strict
smooth subsolution $V_1$ and a strict smooth supersolution $V_2$
to $(\ref{FB1})$ in $\R^n$, such that
\begin{enumerate}
\item $V_1 \leq V_2$ on $\R^n, 0 \in \{V_2>0\} \cap
\{V_1=0\}^\circ$; \item $\p_n V_i >0$ in $\ov{\{V_i>0\}},$ for
$i=1,2$.
\end{enumerate}
Then, for each $R>0$, and $h_R$ sufficiently large, there exists
$u_R$ minimizer of $J$ over $K_R := \{v \in H^1(C_R)| V_1\leq v
\leq V_2 \ \ \text{a.e on $C_R$}, v = V_2 \ \ \text{on $S_R$}\},$
such that
\begin{enumerate}\renewcommand{\theenumi}{\alph{enumi}}\item $u_R$ is
monotone increasing in $C_R^+(u)$ in the $x_{n}$ direction;
\item $u_R$ is Lipschitz continuous, non-degenerate, $(I)$ non-degenerate, and
satisfies the density property $(D)$ in $C_R;$ \item $u_R$ is a
viscosity solution to $(\ref{FB1})$ in $C_R.$\end{enumerate}
\end{thm}

\

\noindent\textbf{Remark 1.} The height $h_R$ of the cylinder
$C_R$, must be such that the free boundaries of $V_1$ and $V_2$
``exit" from the sides of the cylinder. For example, let $\p \{V_i
>0\} = \{(x',\phi_i(x')), x' \in \R^{n-1} \}.$ Then we can
choose
\[h_R \geq \max\{2d_R(V_1),2d_R(V_2), R\}, \ \ \
d_R(V_i)= \max_{\mathcal{B}_R(0)}{|\phi_i(x')|}, i=1,2.\]

\

\noindent\textbf{Remark 2.}  Under the assumptions of Theorem
\ref{localex}, $V_1 < V_2$ on $\{V_2 > 0\}$, and $F(V_1) \cap
F(V_2) = \emptyset$.

\textsc{Proof.} Assume that there exists $x \in\{V_2
> 0\}$ such that, $V_1(x)=V_2(x)$. Then, since $V_1$ is
subharmonic in $\{V_2>0\}$, the maximum principle implies $V_1
\equiv V_2$, which contradicts the fact that $V_1$ is a strict
subsolution and $V_2$ is a strict supersolution. Analogously,
suppose $x \in F(V_1) \cap F(V_2) $, and let $B \subset \{V_1>0\}$
be a ball tangent to $F(V_1)$ at $x$. Then, by Hopf's lemma,
$\p_{\nu} (V_1-V_2)<0$, with $\nu$ inner normal derivative to $\p
B$ at $x$. Again, this contradicts the fact that $V_1$ is a strict
subsolution, and $V_2$ is a strict supersolution. \qed

\

\noindent\textbf{Remark 3.} One could prove that if $u \geq 0$ is
a Lipschitz continuous function on $\Omega$, and $u$ is harmonic
on $\Omega^+(u),$ then the two notions of non-degeneracy given in
Definition \ref{nond} and Definition \ref{Inond} are equivalent
(see also Lemma 7 in \cite{C3}). We do not prove this fact here,
since according to part (b) in Theorem \ref{localex}, we can prove
directly that $u_R$ satisfies both properties. Therefore, we think
of these two definitions as two versions of the same property.

\

Let $W_R := C^+_R(V_2) \cap C_R^-(V_1)$ be the region between the
two smooth (non-touching) graphs $F(V_1), F(V_2)$. Observe that if
$u_R$ minimizes $J$ over $K_R$, then it minimizes $J(\cdot, W_R)$
among all competitors which equal $u_R$ on $\partial W_R$. Hence,
$u_R$ is Lipschitz continuous, non-degenerate and satisfies the
density property $(D)$ in $W_R$  (see \cite{AC}). Thus, in the
proof of part (b) in Theorem \ref{localex}, we only need to
analyze what happens when $F(u_R)$ is ``close" to the fixed
boundary of $W_R$, that is either to $F(V_1)$ or to $F(V_2).$
Whenever we are ``away" from either $F(V_1)$ or $F(V_2)$ the
desired properties will follows with the same arguments as in
\cite{AC} or \cite{ACF}, while near $F(V_1)$ or $F(V_2)$, the free
boundary $F(u_R)$ will inherit the same good properties of the
smooth graphs $F(V_1), F(V_2).$

Moreover, according to Lemma 7 in \cite{C2}, $u_R$ is also a
viscosity solution to \eqref{FB1} in $W_R.$ Thus, in order to show
part (c) in Theorem \ref{localex} it is enough to show that in
fact $F(u_R) \subset W_R$, that is $F(u_R)$ does not touch neither
$F(V_1)$ nor $F(V_2)$.

\

\textsc{Proof of Theorem \ref{localex}.} We divide the proof in
four steps.

\

\noindent\textbf{\textit{Step 1: Existence of monotone
minimizers.}} \comment{We recall the notion of monotone
rearrangement. Let $D \subset \R^n$ be a compact set. For each $x'
\in \R^{n-1}$ we introduce the notation
\[D(x')= D \cap \{(x',x_n)| x_n \in \R\}.\]
Assume that $D$ is convex in $x_n$, i.e. for each $x' \in
\R^{n-1}$, $D(x')$ is either empty, or consists of a single closed
interval.
 For a given $b \geq 0$, we define
  \[D^*(x'):=
  \begin{cases}
    \{(x',x_n)\in \R^n| -b \leq x_n \leq |D(x')| - b\} & \text{if $D(x')\neq \emptyset$}, \\
    \ & \  \\
    \emptyset & \text{otherwise}
  \end{cases}
\] and \[D^* : = \bigcup_{x' \in D'}D^*(x'),\]where $D' \subset
\R^{n-1}$ is the set of those $x' \in \R^{n-1}$ for which $D(x')$
is non empty. Let $C(a,b)=\mathcal{B}_a(0) \times \{|x_n|<b\}$,
$a,b>0$, be a cylinder in $\R^n$, and let $u$ be a Lebesgue
measurable function on $\ov{C(a,b)}$. We define the monotone
(decreasing) rearrangement $u^*$ of $u$, in the direction $x_n$,
by the following formula:
\[u^*(x):= \sup\{k \in \R | x \in (C(a,b)_k)^*\}\] for $x \in
\ov{C(a,b)}$. Here we are using the notation $G_k:= \{x \in \Omega
| u(x)\geq k \}$, for any function $u$ defined on a measurable
subset $G$ of $\R^n$. The function $u^*$ is monotone decreasing in
the $x_n$-direction, and $u$ and $u^*$ are equimeasurable, that
is, for all $k \in \R$ \[|\{u \geq k\}| = |\{u^* \geq k\}|.\]
Moreover, the mapping $u \rightarrow u^*$ is order preserving,
i.e, $u \leq v$ implies $u^* \leq v^*.$ One can define similarly
the concept of monotone increasing rearrangement, which we will
still denote by $u^*$.}We recall that for any  Lebesgue measurable
function $u$ on $\ov{C_R}$ one can define the monotone
(non-decreasing) rearrangement $u^*$ of $u$, in the direction
$x_n$. The function $u^*$ is monotone non-decreasing in the
$x_n$-direction, and $u$ and $u^*$ are equimeasurable, that is,
for all $k \in \R$
\[|\{u \geq k\}| = |\{u^* \geq k\}|.\] Moreover, the mapping $u
\rightarrow u^*$ is order preserving, i.e, $u \leq v$ implies $u^*
\leq v^*.$ The following proposition holds (for details we refer
the reader to \cite{K}).
\begin{prop}\label{rearrang}Let $u \in W^{1,p}(C_R)$, $1<p<\infty$. Then $u^* \in
W^{1,p}(C_R)$ and we have
\begin{equation}\int_{C_R}{|\nabla u|^p dx} \geq \int_{C_R}{|\nabla u^*|^p dx.}\end{equation}
\end{prop}

Now, consider the energy functional $J_R(\cdot)=J(\cdot,C_R)$
defined in \eqref{J}. Denote by $K_R$ the following closed and
convex subset of $H^1(C_R)$,
\[K_R := \{v \in H^1(C_R)| V_1\leq v \leq V_2  \ \ \text{a.e on $C_R$}, v = V_2 \ \ \text{on $S_R$}\}.\]
 The following existence theorem holds.

\begin{thm}\label{exist}
The minimum of $J_R$ over $K_R$ is achieved at $u_R$, which is a
monotone non-decreasing function in $C_R^+(u)$ in the
$x_n$-direction.
\end{thm}

\textsc{Proof.} Since $J_R$ is non-negative, there exists a
minimizing sequence $u_m$, that is
\[u_m \in K_R, \ J_R(u_m)\rightarrow \alpha \equiv \inf_{v \in K_R}{J_R(v)}, \ 0\leq \alpha \leq J_R(V_2)<\infty.\]
The sequence $\{u_m\}$ is uniformly bounded in $H^1(C_R)$. Indeed,
\[||\nabla u_m||_2^2 \leq J(V_2), \ \ ||u_m||_2 \leq ||V_2||_2 .\]
Therefore, we can extract a subsequence, which we will still
denote by $\{u_m\}$, such that $u_m \rightarrow u \in K_R$, weakly
in $H^1(C_R)$. It is known that $J_R$ is lower semicontinuous (see
\cite{AC}) with respect to weak $H^1$ convergence, that is,
\[\liminf_{m\rightarrow \infty}{J_R(u_m) \geq J_R(u)}.\]
%Indeed,
%\begin{equation}\nonumber \int_{C_R}|\nabla u_m|^2 \geq  \int_{C_R}|\nabla
%u|^2 + 2\int_{C_R}{\nabla (u_m - u)\cdot \nabla u,}\end{equation}
%and the right hand side tends to $0$, for $m\rightarrow \infty.$
%Moreover, for each $\eps > 0$, up to extracting a subsequence,
%\begin{align*}
%u_m &\rightarrow u, \ \ \  \text{a.e. on $C_R$}\\
%u_m &\rightarrow u, \ \ \  \text{uniformly on $(C_R \setminus W)$,
%with $|W|< \epsilon$.}
%\end{align*} Thus, for $m$ large, we have
%\[\int_{C_R}{\chi_{\{u_m >0\}}} \geq \int_{C_R\setminus
%W}{\chi_{\{u > \epsilon\}}} \geq \int_{C_R}{\chi_{\{u
%>\epsilon\}}} - \epsilon,
%\] hence \[\liminf_{m\rightarrow \infty}{\int_{C_R}{\chi_{\{u_m >0\}}}} \geq
%\int_{C_R}{\chi_{\{u>0\}}}.\]
This immediately implies that $u$ is a minimizer for $J_R$ over
$K_R$. Now, let $u^*$ be the monotone rearrangement of $u$. Then,
using Proposition \ref{rearrang}, together with the
equimeasurability of rearrangements, we get that
\[J_R(u^*)\leq J_R(u).\] Moreover, the order preserving property
implies that $u^* \in K_R$. Hence $u^*$ is the desired minimizer,
monotone non-decreasing in the $x_n$ direction. \qed

\

We will henceforth denote by $u_R$, a minimizer of $J_R$ over
$K_R$, which is monotone non-decreasing in the $x_n$ direction.

\comment{

\noindent\textbf{Remark.} In this step, we do not need $V_1,V_2$
to be respectively a strict subsolution and a strict
supersolution.}

\

\noindent\textbf{\textit{Step 2: Continuity and harmonicity of
monotone minimizers.}} We wish to show that $u_R$ is harmonic in
its positive phase. This is achieved via standard techniques.

\begin{lem}\label{harm} $u_R$ is continuous in $C^{+}_R(V_2)$, and harmonic in
$C^{+}_R(u_R).$ In particular, $u_R$ is monotone increasing in
$C^{+}_R(u_R)$ in the $x_n$ direction.
\end{lem}

\textsc{Proof.} Let $D$ be a compact subset of $C_R^+(V_2)$, and
let $B_\rho$ be a ball of radius $\rho$ in $D$. Denote by $v_\rho$
the harmonic replacement of $u_R$ on $B_\rho$, that is the
harmonic function in $B_\rho$ which equals $u_R$ on $\partial
B_\rho$. Assume that $v_\rho$ is extended to be $u_R$ outside
$B_\rho$. Since $0 \leq u_R \leq V_2$ a.e., we have $0 \leq
(v_\rho - V_2)^+ \leq (v_\rho - u_R)^+$. Hence $(v_\rho -V_2)^+
\in H^1_0(B_\rho)$ (see \cite{AH}). Therefore, by the weak maximum
principle (see \cite{GT}) we obtain $v_\rho \leq V_2$ a.e. on
$B_\rho$. Analogously, we get $V_1 \leq v_\rho$ a.e. on $B_\rho$.
Since $u_R$ minimizes $J(\cdot,B_\rho)$ among all competitors $v$,
such that $V_1 \leq v \leq V_2$, and $v=u_R$ on $\partial B_\rho,$
we get that
\[\int_{B_\rho}{(|\nabla u_R|^2 + \chi_{\{u_R > 0\}})} \leq \int_{B_\rho}{(|\nabla v_\rho|^2 + \chi_{\{v_\rho >
0\}})}.\] Therefore, \[\int_{B_\rho}(|\nabla u_R|^2 - |\nabla
v_\rho|^2) \leq K \rho^n.\] A standard iterative argument (see
\cite{ACF}) then implies the desired continuity.

\noindent Now, take $\ov{x} \in C_R^+(u_R)$. By continuity, there
exists $r>0$ such that $B_r(\ov{x}) \subset C_R^+(u_R)$. Let
$w_r$, be the harmonic replacement of $u_R$ on $B_r(\ov{x})$.
Since $w_r$ minimizes the Dirichlet integral and $w_r >0$ on
$B_r(\ov{x})$, we get that
\[J(w_r,B_r(\ov{x})) \leq J(u_R, B_r(\ov{x})).\]As before, the
minimality of $u_R$ implies that the reverse inequality holds as
well. Hence \[\int_{B_r(\ov{x})}{|\nabla w_r|^2} =
\int_{B_r(\ov{x})}{|\nabla u_R|^2}.\] By uniqueness of the
Dirichlet minimizer we obtain then $u_R= w_r$ on
$B_r(\ov{x})$.\qed

\

\noindent \textbf{\textit{Step 3: Lipschitz continuity,
non-degeneracy, and density property of monotone minimizers.}} We
start by proving the Lipschitz continuity of $u_R$ in $C_R$. Set
$d(x,K)=\mathrm{dist}(x,K),$ for any closed set $K$. In
particular,$$d(x)=d(x,F(u_R)).$$

\begin{lem}\label{Lip}$u_R$ is Lipschitz continuous in $C_R$, with universal
constant on each $D \Subset C_R$. In particular, for every $D
\Subset D' \Subset C_R,$ there exists $K>0$ depending on $D, n$
and on the $C^1$ norm  of $V_2$ on $D'$, such that,
\[ u_R(x) \leq K d(x), \ \ \ \text{for all $x \in D$}.\]
\end{lem}

\textsc{Proof.} Let $x_0 \in D \Subset D' \Subset C_R$, with
$u(x_0) >0$, and let $B_r = B_r(x_0)$ be the maximum ball
contained in $D' \cap \{u>0\}$. If $\p B_r$ touches $\p D'$, then
$r \geq \text{dist}(D,D')$, and we can apply interior regularity
together with the fact that $ u_R \leq V_2$, in order to show
$|\nabla u_R|(x_0) \leq K$. Otherwise, $\p B_r$ touches $F(u_R)$
at a point $x_1$.\newline
 We distinguish two cases.

\vspace{2mm}

\noindent \textit{(a)} $d(x_1,F(V_2))>r/2.$ Then $B_{r/2}(x_1)
\subset C_R^+(V_2)$ and we can apply the argument of Lemma 3.2
from \cite{AC} to conclude the following, \comment{ We replace
$u_R$ in $B_{r/2}(x_1)$ by the harmonic function $v$ with boundary
values $u_R$. Then, by the maximum principle, $V_1 \leq v \leq
V_2$ on $B_{r/2}(x_1)$, and also $v
> 0$ in $B_{r/2}(x_1)$. Thus, the minimality of $u_R$ yields
\[\int_{B_{r/2}(x_1)}|\nabla(u_R - v)|^2 \leq \int_{B_{r/2}(x_1)}{\chi_{\{u_R=0\}}}.\]
The right hand side can be estimated as in \cite{AC}, Lemma 3.2.
One gets
\[\left(\int_{B_{r/2}(x_1)}{\chi_{\{u_R=0\}}}\right)(\ov{u_R})^2 \leq K r^2 \int_{B_{r/2}(x_1)}|\nabla(u_R - v)|^2, \]
with $K>0$ dimensional constant and $\ov{u_R}$ the average of
$u_R$ on $B_{r/2}(x_1)$. Combining the two estimates above, and
the fact that $x_1 \in F(u_R)$, we obtain
\[\ov{u_R} \leq K r, \] that is }\[\fint_{B_{r/2}(x_1)}{u_R dx} \leq Kr.\]
Now, let $\ov{x}$ be on the ray from $x_0$ to $x_1$, at distance
$r/4$ from $x_1$. Then, by Harnack inequality, and the mean value
property for $u_R$, we get (with $K$ changing for each inequality)
\[u_R(x_0) \leq K u_R(\ov{x})= K \fint_{B_{r/4}(\ov{x})}{u_R} \leq
K \fint_{B_{r/2}(x_1)}{u_R} \leq K r. \]

\vspace{2mm}

\noindent \textit{(b)} $d(x_1,F(V_2))\leq r/2.$ Then,
\[u_R(x_0) \leq V_2(x_0) \leq K d(x_0,F(V_2)) \leq K {|x_0-x_1|+ d(x_1,F(V_2))} \leq K
r,\]where $K$ depends on the Lipschitz norm of $V_2$ on $D'$.

\noindent Now, denote by $v(x)=u_R(r x + x_0)/r$. Then, $\Delta v
= 0$ and by Harnack inequality, $v(x) \leq K$ on $B_{1/2}(0)$. By
interior regularity, $|\nabla u| \leq K' $ on $B_{1/4}(0)$, with
$K'$ dimensional constant. Rescaling back to $u_R$, we obtain
$|\nabla u_R| \leq K'$, on $B_{r/4}(x_0)$, which implies the
desired Lipschitz continuity. \qed

\

\begin{cor}\label{exp1}$u_R$ is a Lipschitz continuous subharmonic function in $C_R.$
\end{cor}

\comment{The following result can be found in \cite{C2}.
\begin{lem}\label{exp} Let  $\Omega_1$ (resp. $\Omega_2$) be such
that \[\Omega_1 \cap B_1(0) \supset \{x_n > 0\} \cap B_1(0),
(\text{resp.} \ \Omega_2 \cap B_1(0) \subset \{x_n > 0\} \cap
B_1(0)).\] Assume that $u$ is a Lipschitz positive harmonic
function in $\Omega_1$ (resp. $\Omega_2$) vanishing on $\partial
\Omega_1$ (resp. $\partial \Omega_2$) and assume that \[\ov{B_1}
\cap
\partial \Omega_i \cap \{x_n=0\}=\{0\}.\]
Then, near zero, $u$ has the asymptotic development \[u(x)= \alpha
x_n^+ + o(|x|) \ \ \text{on $\{x_n>0\}$,}\]with $\alpha \geq 0$.
Furthermore, $\alpha > 0$ for $\Omega_1.$
\end{lem}

Lemmas \ref{harm}, \ref{Lip}, and Lemma \ref{exp}, imply that,
near a regular free boundary point, $u_R$ has the desired
expansion, as in (ii) Definition \ref{defsol}, with $\alpha \geq
0.$ In particular, $\alpha > 0$ at points where an exterior ball
condition is satisfied. We now wish to prove that $\alpha=1.$
Towards this aim it suffices to show that $F(u_R)$ does not touch
$F(V_1)$ and $F(V_2)$. Then, we can apply the argument in Theorem
... \cite{C2}, to conclude that $\alpha=1.$}

We now prove a non-degeneracy result. \comment{We will need the
following auxiliary lemma from \cite{C3}.

\begin{lem}\label{degV1} Let $v$ be a Lipschitz non-degenerate function in $\ov{\Omega} \cap
B_1(0),$ satisfying $\Delta v \geq 0, v=0$ on $\p \Omega \cap
B_1(0).$ Assume further that $0 \in \p \Omega.$ Let $v(x_0) \geq C
d(x_0, \p \Omega),$ for $x_0 \in B_{1/2}(0),$ then, for $\rho \leq
1/4,$ we have
\[ \sup_{B_\rho{(0)}}{v} \geq C \rho.\]
\end{lem}

We are now ready to prove the non-degeneracy property of $u_R.$}

\begin{lem}\label{nondeg} $u_R$  is non-degenerate on
$C_R$, i.e., for every $D \Subset D' \Subset C_R$, there exists
$\ov{K}>0$ depending on $D, n$ and on the $C^2$ norm of $V_1$ on
$D'$, such that\[\ov{K} d(x)\leq u_R(x),\] for all $x \in
D^+(u_R)$, such that $B_{d(x)}(x) \subset D.$
\end{lem}

\textsc{Proof.} Let $x_0 \in D^+(u_R)$, and denote by $r=d(x_0)$.
Assume $B_r(x_0) \subset D.$ We distinguish two cases.

\vspace{2mm}

 \noindent\textit{(a)} $d(x_0, F(V_1)) > r/2.$ If
$x_0 \in C_R^-(V_1),$ then $B_{r/2}(x_0) \subset C_R^-(V_1)$ and
we can proceed as in \cite{ACF} Theorem 3.1, to conclude that

\[u_R(x_0)= \fint_{B_{r/2}(x_0)}{u_R} \geq \ov{K} r. \]

\noindent If $V_1(x_0)>0,$ then
\[u_R(x_0) \geq V_1(x_0) \geq K d(x_0,F(V_1)) \geq \ov{K} r,\]
with $K$ depending on the $C^2$ norm of $V_1$ on $D'$.

\vspace{2mm}

 \noindent\textit{(b)} $d(x_0, F(V_1))=|x_0-x_1| \leq r/2.$ Then $B_{r/2}(x_1) \subset
B_r(x_0)$, hence by Harnack inequality,
\[u_R(x_0) \geq K \sup_{B_{r/4}(x_1)}u_R \geq K\sup_{B_{r/4}(x_1)}V_1 \geq
\ov{K} r,\]again with $\ov{K}$ depending on the $C^2$ norm of
$V_1$ on $D'$. \qed

\

We wish to prove a density property for free boundary points.
Towards this aim, we will need to reformulate our non-degeneracy
property in the following way:
\begin{cor}\label{intnondeg}$u_R$ is $(I)$ non-degenerate on
$C_R$, i.e. for any $D \Subset C_R,$ there exist a constant $K,$
such that, for any ball $B_r \subset D$ centered at a free
boundary point,\[\fint_{B_r}{u_R} \geq K r.\]
\end{cor}
The corollary above can be deduced by the arguments in the proof
of Lemma \ref{nondeg}. We are now ready to derive the following
density property.

\begin{lem}\label{dens}$u_R$ satisfies the density property $(D)$ on $C_R$, i.e.
for any $G \Subset G' \Subset C_R,$ there exist a constant $c<1,$
depending on $G, n$ and on the Lipschitz constant of $F(V_2)$ on
$G'$, such that for any ball $B_r \subset G$ centered at a free
boundary point,
\[c \leq \frac{|B_r \cap \{u_R>0\}|}{|B_r|} \leq 1-c.\]\end{lem}

 \textsc{Proof.} Assume $B_r$ is centered at $0.$
By Corollary \ref{intnondeg}, there exists $y \in \p B_{r/2}$ such
that, $u(y) \geq Kr/2.$ By Lipschitz continuity, for any $z \in
B_{kr}(y)$ we have:
\[u_R(z) \geq u_R(y) - C|z-y|> Kr/2 - Ckr >0\]
as long as $k$ is sufficiently small. Hence $B_{kr} (y) \subset
B_r \cap \{u_R>0\}$, from which the desired lower bound follows.
In order to get the upper bound, we distinguish two cases.

\vspace{2mm}

\noindent \textit{(a)} $d(0,F(V_2))=|x_0|\leq r/2.$ Then
$B_{r/2}(x_0) \subset B_r.$ Hence,\[|\{u=0\}^\circ \cap B_r| \geq
|\{V_2=0\}^\circ \cap B_{r/2}(x_0)|\geq C |B_{r/2}(x_0)| \geq C'
|B_r|.
\]

\vspace{2mm}

\noindent \textit{(b)} $d(0,F(V_2))=|x_0|> r/2.$ Then, $B_{r/2}(0)
\subset \{V_2>0\}$. Hence we can replace $u_R$ with its harmonic
replacement on $B_{r/2}(0),$ and proceed as in \cite{AC}, Lemma
3.7. \qed

\

\noindent \textbf{Remark.} One could prove the statements of Lemma
\ref{Lip} and Lemma \ref{nondeg} with the constants $K, \ov{K}$
independent of $V_2$, and $V_1$ respectively.

\

Lemmas \ref{Lip}-\ref{nondeg}-\ref{dens} and Corollary
\ref{intnondeg} prove the statement in part (b) of our Theorem.

\

\noindent\textbf{\textit{Step 4: Non-touching of the free
boundaries.}} We are now ready to prove the following statement,
using blowing-up techniques.

\begin{lem}\label{FBnotouch1}$F(u_R)$ does not intersect neither
$F(V_1)$ nor $F(V_2).$
\end{lem}

We start by proving that $u_R$ is separated from $V_1$ and $V_2$
in its positive phase.

\begin{lem}\label{strictineq} $u_R < V_2$ in $C_{R}^+(V_2),$ and $V_1 <
u_R$ in $C_{R}^+(u_R).$
\end{lem}
\textsc{Proof.} Assume $u_R(x)= V_2(x)$ at some point $x \in
C_R^+(V_2)$, then the strong maximum principle implies that $u_R
\equiv V_2$ on $C_R^+(V_2)$, hence $u_R \equiv V_2$ on $C_R$. We
want to show that this contradicts the fact that $u_R$ minimizes
$J_R$ on $K_R.$ Let $g \in C^{\infty}_0(C_R), g \leq 0$. For
$\epsilon > 0,$ set $y_{\epsilon}(x)= x+\epsilon g e_n$ and
$V_{\epsilon}(x)=u_R(y_{\epsilon}(x)).$ For $\epsilon$
sufficiently small, the monotonicity of $V_2$ in the
$x_n$-direction and the fact that $V_1 < V_2$ in the positive
phase of $V_2$, imply that $V_{\epsilon} \in K_R$. Therefore,
using that Det$(y_{\epsilon}(x))= 1+\epsilon \nabla \cdot ge_n
+o(\epsilon^2)$, we get that
\begin{multline}\nonumber 0 \leq J_R(V_{\epsilon})- J_R(u_R)=\\=
\epsilon \left\{\int_{C_R}{-(|\nabla u_R|^2 +
\chi_{\{u_R>0\}})\nabla \cdot ge_n +(2\nabla u_R D ge_n \nabla
u_R)}\right\} + o(\epsilon^2).
\end{multline}
Therefore using Lemma \ref{harm}, we obtain ($u_R \equiv V_2$ and
$V_2$ smooth)
\begin{multline}\nonumber 0 \geq \int_{\{V_2>0\}}{\nabla
\cdot ((|\nabla V_2|^2 + 1)ge_n - 2g e_n\cdot \nabla V_2 \nabla
V_2)}=\\= - \int_{\partial\{V_2>0\}}{((|\nabla V_2|^2 + 1)ge_n -
2g e_n\cdot \nabla V_2 \nabla V_2)\cdot \nu}=\\= -
\int_{\partial\{V_2>0\}}{(1-|\nabla V_2|^2)g \nu_n}\end{multline}
for all function $g$ as above, and $\nu$ the inner unit normal to
$\p \{V_2>0\}$. This contradicts the strict supersolution property
of $V_2$.

Assuming now, $u_R(x)=V_1(x)$ at some point $x \in C^{+}_R(u_R)$,
then the contradiction follows immediately by the fact that $V_1 <
V_2$ on $\{V_2 > 0\}$, and $u_R = V_2$ on $S_R$.\qed

\

\textsc{Proof of Lemma \ref{FBnotouch1}.} We prove that $F(u_R)$
does not touch $F(V_2)$. The proof that $F(u_R)$ does not
intersect $F(V_1)$ follows by similar arguments.

\noindent Assume by contradiction that there exists $x_0 \in
F(u_R)\cap F(V_2).$ \comment{Let $B_r \subset C_R^+(V_2)$ be a
ball tangent at $x_0$ to $F(V_2)$. By Lemma \ref{strictineq} and
Corollary \ref{exp1}, we can apply Hopf's Lemma to the subharmonic
function $v_R= u_R - V_2$ and conclude that
 \begin{equation}\label{hopf}\liminf_{x \rightarrow x_0}
(v_R(x)/|x-x_0|) < 0.\end{equation} }Then, $F(u_R)$ has a tangent
ball from the zero side at $x_0$, and according to Lemma A1 in
\cite{C2}, $u_R(x) = a (x-x_0, \nu)^+ + o(|x-x_0|)$, near $x_0$,
from the positive side of $u_R$, with $\nu$ the inner normal to
$F(V_2)$ at $x_0$. Furthermore, by non-degeneracy (Lemma
\ref{nondeg}), $a>0$. Let $B_{\rho_k}(x_0)$ be a sequence of balls
with $ \rho_k \rightarrow 0$ such that $u_k(x):=
\frac{1}{\rho_k}u_R(x_0+\rho_k x)$ blows up to $U(x),$ and
$V_k(x):= \frac{1}{\rho_k}V_2(x_0+\rho_k x)$ blows up to
$\ov{V_2}.$ Thus, on the unit ball $B$, $U(x)= a (x,\nu)^+$ and
$\ov{V_2}(x)= b (x,\nu)^+$ , and
\begin{equation}\label{ab}0<a \leq b<1.\end{equation}

Let $B_{\delta}(x_0)$ be a ball centered at $x_0$, with $\delta$
small enough so that $B_{\delta}(x_0) \subset C_R^-(V_1).$ Let $g
\in C^{\infty}_0(B_{\delta}(x_0)), g \leq 0$. For $\epsilon
> 0,$ set $y_{\epsilon}(x)= x+\epsilon g e_n$ and
$u_{\epsilon}(x)=u_R(y_{\epsilon}(x)).$ For $\epsilon$
sufficiently small, the monotonicity of $u_R$ in the
$x_n$-direction and the fact that $B_{\delta}(x_0) \subset
C_R^-(V_1)$, imply that $V_1 \leq u_{\epsilon} \leq V_2$ on
$B_{\delta}(x_0)$, and hence
$$0 \leq J(u_{\epsilon}, B_{\delta}(x_0))- J(u_R,
B_{\delta}(x_0)).$$ After rescaling appropriately and passing to
the blow-up limit (using Lemma \ref{blowup}), the inequality above
implies that
$$0 \leq J(U_{\epsilon}, B)- J(U,
B),$$ where $U_\epsilon (x) = U (x+ \epsilon \tilde{g}e_n)$, for
any $\tilde{g}$ compactly supported in the unit ball $B$ and
$\tilde{g} \leq 0.$ The same computations as in Lemma
\ref{strictineq} give\begin{equation} 0 \leq \lim_{\epsilon
\rightarrow 0} \frac{1}{\epsilon}(J(U_{\epsilon}, B)- J(U, B)) =
\int_{\partial\{U>0\}}{(1-a^2)\tilde{g} \nu_n}\end{equation} for
all function $\tilde{g}$ as above. Thus, we contradict \eqref{ab}.

\comment{Let us prove that $U$ is an absolute minimum for
$J(\cdot,B)$, among all competitors $v \leq \ov{V_2}, v <
\ov{V_2}$ in $\{\ov{V_2}>0\},v =U$ on $\partial B$. Let $v$ be an
admissible competitor and define
\[v_k=v+(1-\eta)(u_k -U)\]for $\eta \in C^{\infty}_{0}(B), 0 \leq \eta \leq
1$. Set $w_k = v_k^+$. Then, $w_k = u_k$ on $\partial B,$ and for
$k$ large enough $w_k(x) \geq 0 = \frac{1}{\rho_k}V_1(x_0 + \rho_k
x).$ Moreover, using that $v_k - V_k$ converges uniformly to
$v-\ov{V_2}$, and $u_k \leq V_k(x)$, we get that $w_k \leq V_k$,
for $k$ large enough. Therefore $J(u_k,B) \leq J(w_k,B)$ from
which we obtain
\begin{multline}\nonumber \int_{B}{\nabla((U-v)+\eta(u_k -U))\cdot
\nabla((U+v)+(2-\eta)(u_k - U))} +\\+ \int_{B}{(\chi_{\{u_k>0\}} -
\chi_{\{v_k>0\}})} \leq 0.\end{multline} Observe that the
following inequality holds
\[\chi_{\{w_k>0\}} \leq \chi_{\{v>0\}} + \chi_{\{\eta <1\}}\]
Letting $k \rightarrow \infty$, we get $J(U,B) \leq J(v,B) +
|\{\eta <1\}|$, and an appropriate choice of the function $\eta$
gives the desired minimality. Now let $g \in C_{0}^{\infty}(B), g
\geq 0$ and for $\epsilon >0$, set $U_{\epsilon}(x)= U(x-\epsilon
g\nu).$ Applying the same domain variation technique as in Lemma
\ref{strictineq}, we therefore get $a \geq 1$, which is a
contradiction.}\qed

Combining the two Lemmas above, we obtain the following Corollary.

\begin{cor}\label{variat}$u_R$ is a variational solution to $(\ref{FB1})$ in $C_R$.
\end{cor}
\comment{\textsc{Proof.} Let $\eta \in C_0^{\infty}(C_R, \R^n)$,
and $\eps$ small. Define $u_\eps(x)= u_R(\tau_\eps(x)),$ where
$\tau_\eps(x) = x + \eps \eta(x)$. Then, Lemmas \ref{FBnotouch1}
and \ref{strictineq}, guarantee that $u_\eps \in K_R.$ By the same
computations as in Lemma \ref{strictineq}, we therefore get the
desired limiting equality. \qed}

\

%Using the corollary above, we can now prove the following result.

%\begin{lem}\label{FBnotouch2}$F(u_R)$ does not intersect $F(V_1).$
%\end{lem}
\comment{\textsc{Proof.} Assume by contradiction that there exists
$x_0 \in F(u_R) \cap F(V_1)$, and let denote by $u_0$ a blow-up of
$u_R$ around $x_0$. From Lemma \ref{blowup}, we deduce that we can
pass to the limit in the definition of variational solution, hence
$u_0$ is a variational solution to the one-phase free boundary
problem (\ref{FB1}) on any compact of $\R^n$. Moreover, $u_0$ is
harmonic in its positive phase, hence as in Corollary
\ref{variat}, $u_0$ satisfies the equality:
\begin{equation}\label{hadamard1}\lim_{\eps \rightarrow 0}
\int_{\p \{u_0>\eps\}} {(|\nabla u_0|^2 -1)\eta \cdot \nu =0}.\end{equation}

\noindent Since $x_0 \in F(u_R) \cap F(V_1),$ $u_R$ has an
asymptotic expansion around $x_0$, $u_R(x) = a (x-x_0, \nu)^+ +
o(|x-x_0|)$, with $a > 0$, and $\nu$ the inner unit normal to
$F(V_1)$ at $x_0$. Thus, applying formula (\ref{hadamard1}) to the
blow-up limit $u_0(x)=a (x,\nu)^+ $, we get $a=1$. Since $V_1$ is
a strict subsolution, Hopf's lemma implies $a>1$. We have reached
a contradiction, hence $F(u_R)$ and $F(V_1)$ cannot touch. \qed}

\section{Local regularity theory}

Throughout this section we will denote by $u_R$ a viscosity
solution to the problem
\begin{equation}
  \begin{cases}\label{FB2}
    \Delta u = 0 & \text{in $\Omega^+(u)$}, \\
   |\nabla u|  = 1 & \text{on $F(u)$},
  \end{cases}
\end{equation}
in the cylinder $C_R = \mathcal{B}_R(0) \times \{|x_n| < h_R\}$.
Moreover, $u_R$ is monotone increasing in $C_R^+(u_R)$ in the
$x_n$ direction, $V_1 \leq u_R \leq V_2$ in $C_R$, and $u_R=V_2$
on the sides $S_R= \p \mathcal{B}_R(0)\times \{|x_n| \leq h_R\}$.
Recall that $V_1$ is
 a strict
smooth subsolution and $V_2$ is a strict smooth supersolution  to
(\ref{FB2}) in $\R^{n}$, such that
\begin{enumerate}
\item $V_1 \leq V_2$ on $\R^{n}, 0 \in \{V_2>0\}\cap \{V_1 =0\}^\circ;$
\item $\p_{n} V_i >0$ in $\ov{\{V_i>0\}},$ for $i=1,2$.
\end{enumerate}

\

The existence of $u_R$ on a sufficiently tall cylinder is
guaranteed by Theorem \ref{localex}. Let $u_R$ be extended to zero
on $\{(x',x_n):|x'|\leq R, x_n \leq - h_R\}.$

 In this section we wish to prove the following regularity property for the
 free boundary of $u_R.$

\begin{thm}\label{regvis} For $h_R$ sufficiently large,
$F(u_R)$ is a Lipschitz graph in the $x_n$ direction.
\end{thm}

\noindent\textbf{Remark 1.} The results in \cite{C1} and
\cite{KN}, imply
 that
$F(u_R)$ is smooth in the interior.

\vspace{2mm}

\noindent\textbf{Remark 2.} The Lipschitz constant of $F(u_R)$ on
$C_R$ depends on the $C^2$-norm of $V_2$ on $\ov{C_R},$ hence the
result of Theorem \ref{regvis} does not provide a ``local" bound.
Thus, we cannot control the Lipschitz constant on $F(u_R)$ in the
limit as $R \rightarrow \infty$, using this bound. However, the
Lipschitz constant is the same for all viscosity solutions $u_R$
as in the beginning of this section.

\

We start by introducing a particular family of viscosity
subsolutions (see \cite{C1}) to the free boundary problem
\eqref{FB2}.

\begin{lem}\label{conv}
Let $u$ be a viscosity solution to $(\ref{FB2})$ in $\Omega$. Let
$v_t(x) = \sup_{B_t(x)} u(y), t > 0$. Then $v_t$ is a subsolution
to $(\ref{FB2})$ in its domain of definition. Furthermore, any
point of $F(v_t)$ is regular from the positive side.
\end{lem}

We will also need the following results from \cite{C1}.
\begin{lem}\label{notouch}
Let $v \leq u$ be two continuous functions in $\Omega$, $v < u$ in
$\Omega^+(v)$, $v$ a subsolution and $u$ a solution. Let $x_0 \in
F(v) \cap F(u)$. Then $x_0$ cannot be a regular point for $F(v)$
from the positive side.
\end{lem}

\begin{lem}\label{complemma}(Comparison lemma) Let $v_t$ be a family of subsolutions to \eqref{FB2}, continuous in
$\ov{\Omega} \times [0,T]$. Let $u$ be a supersolution to
\eqref{FB2} continuous in $\ov{\Omega}.$ Assume that
\begin{enumerate}
\item $v_0 \leq u$ in $\Omega,$
\item $v_t \leq u$ on $\partial \Omega$ and $v_t < u$ in $\ov{\Omega^+(v_t)} \cap \partial
\Omega,$ for $0 \leq t \leq T,$
\item every point $x_0 \in F(v_t)$ is regular from the positive side,
\item the family $\Omega^+(v_t)$ is continuous.
\end{enumerate}
Then, $v_t \leq u$ in $\Omega$ for every $t \in [0,T].$
\end{lem}

We now proceed to prove the following technical lemma. Here and
henceforth the constant $C$ will depend on the $C^2$-norm of $V_2$
on $\ov{C_R}.$

\begin{lem}\label{barrier}(Existence of a barrier) There exists $\delta$ positive
and small, such that for every $x_0 \in S_R \cap \ov{\{V_2>0\}}
\cap \{|x_n| \leq h'_R\}$, $h'_R < h_R$, there exists a function
$V_{x_0}$ with the following properties
\begin{enumerate}\renewcommand{\theenumi}{\alph{enumi}}
\item $V_{x_0}(x_0)=u_R(x_0)$,
\item $u_R \geq V_{x_0}$
\item for all $x \in S_R \cap \ov{\{V_2>0\}}
\cap \{|x_n| \leq h'_R\}$, in a $\delta$-neighborhood of $x$, the
level set $V_x = u_R(x)$  is a Lipschitz graph in the vertical
direction, with Lipschitz constant bounded by a uniform constant
$C$ independent of $x$.
\end{enumerate}
\end{lem}
\textsc{Proof.} Let $x_0 \in S_R \cap \ov{\{V_2>0\}} \cap \{|x_n|
\leq h'_R\}$. Let $H$ be the vertical half-space tangent to $S_R$
at $x_0$ and containing $C_R$. Without loss of generality we can
assume $H=\{x_1>0\}$. We wish to construct a function $W_{x_0}$
such that $W_{x_0}^+$ is a subsolution and  $V_2 \geq W_{x_0}^+$
on $H$, and also $V_2(x_0)=W_{x_0}(x_0)$. We have
\begin{equation}\label{boundary} V_2(x) \geq V_2(x_0) + l\cdot (x-x_0) - K|x-x_0|^2,\end{equation}
for some positive constant $K$, and vector $l$. Set
$$W_{x_0}(x)= (V_2(x_0) + (l-ae_1)\cdot (x-x_0) -
K|x-x_0|^2)e^{-x_1}.$$  By choosing $a$ sufficiently large, we can
guarantee that $W_{x_0}^+$ is a strict subsolution in $H$, that is
$W_{x_0}$ is subharmonic in its positive phase and $|\nabla
W_{x_0}|>1$ on $F(W_{x_0}).$ Indeed,
\begin{equation}\Delta W_{x_0}= -2e^{-x_1}[l_1-a-2Kx_1+nK] + W_{x_0},\end{equation}
hence on $H \cap \{W_{x_0}>0\}$, if $a> l_1+ nK$ then $\Delta
W_{x_0}>0.$ Also,
\begin{equation}\partial_1 W_{x_0}= (l_1-a-2Kx_1)e^{-x_1} - W_{x_0},\end{equation}
hence on $F(W_{x_0})$, if $a>l_1+1$ then $\partial_1 W_{x_0} <
-1.$

Let $B$ be the (closed) ball on which $W_{x_0}$ is non-negative,
and let $D$ be the $n-1$ dimensional ball  $D:= B \cap
\partial H$. Let $0$ be the center of $D$, and let us consider
rescale of $W_{x_0}$ around $0$, that is:$$W_t(x)= t
W_{x_0}\left(\frac{x}{t}\right).$$ One can easily verify that the
family $W_t^+$ is non-increasing in $t$, in particular $$W_t^+
\leq W_1^+=W_{x_0}^+, \ \ \ \text{for all $0 \leq t \leq 1$}.$$
This follows from the fact that the function $W_{x_0}$ is
decreasing along each ray originating from $x_0$ into the
half-space $H.$ Hence, $V_2 \geq W_t^+$ in $H$ for all $t \in
[0,1].$ Now, since $u=V_2$ on $S_R \subset H$, we can apply lemma
\ref{complemma} to $u$ and the family $W_t^+$, to conclude that $u
\geq W_{x_0}^+$ in $C_R$, which gives the desired claim with
$V_{x_0}=W_{x_0}^+$. \qed

\

 Theorem \ref{regvis} is an immediate corollary to the
following Theorem.

\begin{thm}\label{regvis2}
For $h_R$ sufficiently tall, there exists a small constant
$0<c<1$, such that, for small $s> 0,$
\begin{equation}\label{lip1}
\sup_{B_{c s}(x)} u_R(y - s e_n) \leq u_R(x),
\end{equation}
for all $x \in C_{R-cs}$.
\end{thm}
\textsc{Proof.} We divide the proof in 4 different steps.

\noindent{\it Step 1.} In this step we show that by taking $h_R$
sufficiently tall, then \eqref{lip1} holds away from the free
boundary. Indeed, let $h_R$ be sufficiently large, so that a strip
$T=\{|x'|\leq R, r_1 \leq x_n \leq r_2 \} \Subset C_R^+(V_1),$ and
$F(V_2)$ exits from the sides of $C_R.$ Then,
$$u_R \geq \min_{T} V_1 = C  \ \ \ \text{on T.}$$ Now, let
$(0,r) \in T$ and let $(0, y_n)$ belong to the level set $\{u_R =
C/2\}.$ Then,
\begin{equation}
C/2 \leq u_R(0,r) - u_R(0,y_n)= \int_{y_n}^{r}\partial_{x_n}
u_R(0, t) dt.
\end{equation} Hence $\partial_{x_n} u_R(0,
t) \geq C/r$, for some $t \in (y_n,r).$ Moreover, $u_R$ grows
linearly away from the free boundary, hence a ball centered at
$(0,y_n)$ and of radius comparable to $C$ is contained in the
positive phase of $u_R$. Since $u_R$ is monotone increasing in the
vertical direction, we can cover the segment joining $y_n$ and $r$
with a finite (depending on $R, F(V_1),F(V_2)$) number of balls
with radii comparable to $C$. Thus, since $\partial_n u_R$ is a
positive harmonic function in $C^+_{R}(u_R)$, Harnack's inequality
implies that $\partial_{n} u_R \geq c$ at $(0,r)$ for some
constant $c$.  Now Harnack's inequality up to the boundary implies
that
\begin{equation}\label{top}\p_n u_R \geq  M \ \ \text{on} \ \ T,\end{equation}  for some
constant $M$. In particular, for $r$ fixed between $r_1$ and
$r_2$, there exists a small $c$, such that
\begin{equation}\label{lip2}
\sup_{B_{c s}(x)} u_R(y - s e_n) < u_R(x),
\end{equation}
for all $x \in \{x_n=r\} \cap \{|x'| \leq R-cs\}$, and $s$ small.

\

\noindent{\it Step 2.} In this step we construct a  family of
subsolutions which we wish to compare to the solution $u_R$ via a
continuity argument, in the spirit of Lemma \ref{complemma}. To
guarantee that we can start the argument, we need to introduce the
following quantity, which represents the length of the longest
vertical segment contained in the free boundary of $u_R$:
\[\ov{s} = \sup\{\lambda
> 0 | \ \ \exists \ \ \ov{x}, \ \ \text{s.t.} \ \ (\ov{x} + \nu
e_n) \in F(u_R), \ \ \forall \ \ |\nu|\leq \lambda \}.\] Since
$V_1 \leq u_R \leq V_2$, $0 \leq \ov{s} < +\infty.$ Let $s$ be a
small positive number, and define $u_{s}(x) = u_R(x-(s +
\ov{s})e_n).$ Now, consider the family of subsolutions
$$v_t^{s}(x) = \sup_{B_t(x)} u_{s}(y), \ \ \ t \geq 0 \ \ \ \text{and
small}.$$

\noindent Since $u_R \leq V_2$ in $C_R$, by Lemma \ref{barrier},
we obtain that all the level sets of $u_R$ in a small neighborhood
of the boundary $S^+_R=S_R \cap \ov{\{V_2>0\}} \cap \{x_n \leq r_2
\}$, are trapped between two Lipschitz graphs, with Lip constants
uniformly bounded. Hence, there exists a small constant $c$, such
that if $dist(x, S_R) = cs$, then
\begin{equation}\label{lip3} \sup_{B_{\frac{c}{10} s}(x)} u_R(y -
(s+\ov{s}) e_n) \leq u_R(x).
\end{equation} Call $s'=c s/10.$
Denote by $$\Omega=\{-h_R < x_n < r\} \cap \{x \in C_R : dist(x,
S_R)
> cs\}.$$ Then, according to \eqref{lip2}-\eqref{lip3},
\begin{equation}v_t^s(x) \leq u_R(x) \ \ \ \text{for} \ \ \ x \in
\p \Omega,  \ \ \text{for all} \ \ t \in [0, s'].
\end{equation} Moreover,
\begin{equation}\label{strict1}v_t^s(x) <u_R(x) \
\ \ \text{for} \ \ \ x \in \ov{\Omega^+(v_t^s)} \cap \p \Omega
.\end{equation} Indeed \eqref{lip3} is strict on $
\ov{\Omega^+(v_{s'})} \cap \p \Omega,$ since whenever $u_R(x)$ is
zero, $B_{s'}(x-(s+\ov{s})e_n)$ is contained in the zero phase of
$V_2$, hence it cannot be tangent to $F(u_R)$. Also, the
monotonicity of $u_R$ in the $x_n$ direction, guarantees that
\begin{equation}\label{noempty} v_{0}^{s} \leq u_R \ \ \ \text{on $\Omega.$}\end{equation}
%Moreover, by the definition of $\ov{s},$
%\begin{equation}\label{strictnoempty} v_{0}^{s} < u_R \ \ \ \text{on $\ov{C^+(v_0^s)} \cap \p C.$}\end{equation}
%We will show that, there exists a constant $\tau =
%\tau(s,\ov{s})$, such that $v_{\tau}^{s} < u_R$ on
%$\ov{\Omega^+(v_{\tau}^{s})} \cap \p \Omega$, and $\tau(s,\ov{s})
%\rightarrow \ov{\tau} > 0$, as $s \rightarrow 0$. First, observe
%that, by the definition of $r_1$,
%\[\ov{\Omega^+(v_{t}^{s})} \cap \p \Omega =  \{u_R=c\} \cup (
%\ov{\Omega^+({v_{t}^{s}}}) \cap S_R). \]According to \eqref{top},
%there exists a $\theta_1$ depending on $\inf_{\ov{D_R}} \p_n u_R$,
%such that $v_{\tau_1}^{s} < u_R$ on $\{u=c\}$, for $\tau_1 = (s +
%\ov{s}) \sin \theta_1$. Furthermore, there exists $\theta_2$,
%depending on $V_2, R$, such that
%\[\sup_{B_{\tau_2}(x)} V_2(y - (s + \ov{s})e_n) \leq V_2
%(x)=u_R(x), \ \ \text{for all $x \in S_R$, with$\tau_2 = (s +
%\ov{s}) \sin \theta_2$},
%\] and the inequality is strict on $S_R \cap \{V_2>0\}$. This implies, $v_{\tau_2}^{s} < u_R$ on $(
%\ov{\Omega^+({v_{\tau}^{s}}}) \cap S_R).$ Finally, $\tau =
%\min\{\tau_1, \tau_2 \}.$\vspace{2mm}

\noindent \textit{Step 3.} In this step, we compare $u_R$ with the
family $v_t^s$, using a standard continuity argument. Define
$A=\{t \in [0, s'] | v_t^s \leq u_R, \text{in $\Omega$}\}.$ By
(\ref{noempty}) $A \neq \emptyset$, and by the continuity in $t$
of the family $v_{t}^{s}$, A is closed. We want to prove that $A$
is open, hence $A=[0,s']$.

\noindent Let $t_0 \in A$, then $v_{t_0}^s \leq u_R$ in $\Omega$,
and by \eqref{strict1} and the strong maximum principle we get
$v_{t_0}^s < u_R$ in $\Omega^+(v_{t_0}^{s}).$ If $t_0
> 0$, then by Lemma \ref{conv}, every point of
$F(v_{t_0}^{s})$ is regular from the positive side, and Lemma
\ref{notouch} implies that $F(v_{t_0}^{s}) \cap F(u_R) =
\emptyset$. Hence,
\[\ov{\Omega^+(v_{t_0}^{s})} \subset \{x \in \ov{\Omega} \
\ | \ \ u_R(x) > 0\}.\] The inclusion above, for the case $t_0 =
0$ follows from the definition of $\ov{s}$. By the continuity in
$t$, for $t$ close to $t_0$,
\[\ov{\Omega^+(v_{t}^{s})} \subset \{x \in \ov{\Omega} \ \ |
\ \ u_R(x) > 0\}.\] Thus, $v_{t}^{s} < u_R$ on $\p
\Omega^+(v_{t}^{s})$. Since $v_{t}^{s} - u_R$ achieves its maximum
on the boundary, we then get $v_{t}^{s} < u_R$ on $
\Omega^+(v_{t}^{s})$, from which we conclude that $t \in A$.

\vspace{2mm}

\noindent \textit{Step 4.} In this step we prove the desired
statement \eqref{lip1}. From Step 3 we
have:\begin{equation}\label{lip}\sup_{B_{c\frac{s}{10}} (x)} u_R(y
- (s + \ov{s})e_n) \leq u_R(x), \ \ \text{in
$\Omega.$}\end{equation} Introduce the family of subsolutions
$$w_t(x)= \sup_{B_{c\frac{s}{10}} (x)} u_R(y
- se_n + (t-1)\ov{s}e_n).$$ Similar arguments as in Step 2,
guarantee that hypothesis (ii) in Lemma \ref{complemma} is
satisfied. Moreover, according to \eqref{lip}, $w_t \leq u_R$ for
$t=0$, hence we can apply Lemma \ref{complemma} to conclude that
this inequality is true for $t=1$ as well, which is the desired
claim.
 \qed

\section{Local NTA property.}

Let $u_R$ be as in Theorem \ref{localex}. In this section we prove
the following property.

\begin{thm}\label{nta} $F(u_R)$ is locally NTA.\end{thm}

The proof of this result follows the line of \cite{ACS}, where NTA
regularity for the free boundary of an optimization problem in
heat conduction is proved.

We need the following monotonicity formula from \cite{ACF}.

\begin{thm}\label{monot}
Let $v$ be a continuous function defined on $B=B_R(x_0).$ Suppose
that $v$ is harmonic in the open set $\{x \in B | v(x) \neq 0\}.$
Let $A_1$ and $A_2$ be two different components in $B$ of the set
$\{x \in B | v(x) \neq 0\}$. Assume that for some constant $c >
0$, and any $r, 0 < r <R$, \[|B_r(x_0) \setminus (A_1 \cup A_2)|
\geq c |B_r(x_0)|.\] Define, for $0 <r < R$,\[\phi(r) =
\left(\frac{1}{r^2} \int_{B_r(x_0) \cap A_1} |\nabla v |^2
\rho^{2-n} dx\right)\left(\frac{1}{r^2} \int_{B_r(x_0) \cap A_2}
|\nabla v |^2 \rho^{2-n} dx\right)\] where $\rho = \rho (x) = |x -
x_0|.$ Then, for some positive $\beta$ depending only on the
dimension and the constant $c$, $r^{-\beta}\phi(r)$ is a
non-decreasing function of $r$.
\end{thm}

We also need the following result for harmonic functions, which
can be obtained with an iterative argument in the spirit of Lemma
7 in \cite{C3}. First we introduce a notation. For any real-valued
function $u$ defined on a domain $\Omega \subset \R^n,$ and any $d
\in \R,$ we denote by
$$\Omega^d(u)=\{x \in \Omega | u(x)>d\}, \ \ \ F^d(u)= \{x \in \Omega | u(x)\leq d\}.$$
\begin{lem}\label{acs42} Let $u$ be a non-negative function in
$\Omega$, $u$ harmonic in $\Omega^+(u)$,  Lipschitz continuous and
non-degenerate. Then, for any compact $D \subset \Omega$, there
exist constants $\beta, \gamma>0$ such that, whenever $\ov{x} \in
F(u) \cap D^\circ, B_{2R}(\ov{x}) \subset D, x\in B_{R/2}(\ov{x})
\cap \{u>0\}$ and $A$ is a connected component of
$\Omega^{\frac{u(x)}{2}}(u) \cap B_R(\ov{x})$ containing $x$, then
there exists $B=B(\ov{y}, \gamma R) \subset A,$ such that
$$\fint_{B}|\nabla u|^2 \geq \beta.$$

\end{lem}

With the same notation as in Lemma \ref{acs42}, we have
\begin{cor}\label{lowbound}There exists a constant $\tau$, such that \[\int_A |\nabla u |^2\rho^{2-n} dy > \tau R^2\] where
$\rho=\rho(y)=|y-\ov{x}|.$

\end{cor}

\textsc{Proof of Lemma \ref{acs42}.} Let $x=x_1$ and denote by
$d_1=|x_1-x_0|$ the distance of $x_1$ from
$F^{\frac{u(x_1)}{2}}(u).$ By the Lipschitz continuity and
non-degeneracy of $u$, we have that the quantities,
$$d_1, u(x_1), d(x_1, F(u))$$ are comparable, with constants
depending only on the Lipschitz and non-degeneracy constants of
$u$ on $D.$ We wish to prove that there exists a positive constant
$\delta$ and a point $x_2 \in B(x_1,d_1)$ such that
\begin{equation}\label{ux2}u(x_2)
\geq (1+\delta)u(x_1).\end{equation} Indeed, by the Lipschitz
continuity of $u$, there exists a constant $c$ such that $$u(y)
\leq \frac{3u(x_1)}{4} \ \ \ \ \text{on $B(x_0, cd_1)$}.$$ Hence,
the function $$v := u - \frac{u(x_1)}{2}$$ is a non-negative
harmonic function in $B(x_1,d_1),$ such that $v(x_1)=u(x_1)/2$ and
$v \leq u(x_1)/4$ in $B(x_0,cd_1)$. Since, $$v(x_1) =
\fint_{B_r(x_1)} v, \ \ \ \ r= d_1(1-c/2),$$ we conclude that
there exists $x_2 \in B_r(x_1)$ such that $v(x_2) \geq
(1+2\delta)v(x_1)$ for a fixed positive constant $\delta.$ This,
according to the definition of $v$, implies the desired claim
\eqref{ux2}. Moreover, from \eqref{ux2} we deduce that there
exists $y_1$ on the line joining $x_1$ and $x_2$, such that
$|\nabla u|^2(y_1) > \beta$ for some small constant
$\beta=\beta(\delta).$ Thus, since $|\nabla u|^2$ is subharmonic
in $B(x_1,d_1)$, we get that
\begin{equation}\fint_{B_{\delta_1}(y_1)} |\nabla u|^2 \geq \beta,\end{equation}
for $\delta_1=c d_1/2$ comparable to $d_1,$ and $B(y_1,\delta_1)
\subset B(x_1,d_1.)$

Now, we iterate this argument $n$- times, to find points $x_j$
such that the following quantities are comparable,
\begin{equation}\label{ud}d_j, u(x_j), d(x_j, F(u)),\end{equation} and
\begin{equation}\label{uxj}u(x_{j+1}) \geq (1+\delta) u(x_{j}), \
\ \ \text{for all $j=1,...,n$}.\end{equation} Also,
\begin{equation}\fint_{B_{\delta_j}(y_j)} |\nabla u|^2 \geq c, \ \ \ B_{\delta_j}(y_j)
\subset B(x_j,d_j), j=1,...,n+1\end{equation} with $\delta_j$
comparable to $d_j.$ The number $n$ is determined so that
$B(x_{n+1},d_{n+1})$ is the first ball to exit $B(x_0, R),$ but it
is still included in $B(x_0, 2R).$ Thus, our lemma is proved with
$\ov{y}= y_{n+1},$ if we show that $d_{n+1}$ is comparable to $R$.
From \eqref{uxj} and the fact that the quantities in \eqref{ud}
are comparable, we get that
$$cR \leq c |x_{n+1} - x_1| \leq c\sum_{i=1}^{n}d_i \leq \sum_{i=1}^{n} u(x_i) \leq C u(x_{n+1}) \leq \tilde{C}
d_{n+1}.$$ This concludes our proof.

 \qed

Finally, the proof of Theorem \ref{nta}, is obtained combining the
non-degeneracy and density property of $u_R$, together with the
Harnack chain property from the next Lemma.

\begin{lem}\label{HC}
Let $u$ be a viscosity solution to $(\ref{FB1})$ in $B_1$, such
that $u$ is Lipschitz continuous and nondegenerate, and $u$
satisfies the density property (D). Then, there exists constants
$M,\ov{\delta}$ such that, for any $\epsilon>0$, and for any
$x_1,x_2 \in B_{1/2},$ such that $B(x,\epsilon) \subset
B_{3/4}^+(u)$ and $|x_1-x_2| \leq C\epsilon \leq \ov{\delta}$,
there exist $y_1=x_1,...y_l=x_2$, such that
\begin{enumerate}\renewcommand{\theenumi}{\alph{enumi}}
\item $B_i=B(y_i, \epsilon/M) \subset B_{3/4}^+(u),$ $i=1,...,l$
\item $B_i \cap B_{i+1}\neq\emptyset$, $i=1,...l-1$
\item $l$ independent of $\epsilon, x_1,x_2.$
\end{enumerate}
\end{lem}
\noindent\textsc{Proof.} Assume that, without loss of generality,
\[\tilde{\delta} = \max\{d(x_1,\p B_{1/2}^+(u)), d(x_2,\p
B_{1/2}^+(u)) \}= d(x_2,\p B_{1/2}^+(u)).\]

\noindent If $\tilde{\delta} \geq 2 C \epsilon$, then, $x_1 \in
B(x_2, C \epsilon) \subset B_{1/2}^+(u),$ and we can easily find
the required chain.

\noindent Assume then, $\tilde{\delta} < 2 C \epsilon$ and let
$x_0 \in \p B_{1/2}^+(u)$ be such that $\tilde{\delta}=|x_2-x_0|.$
Set $r_0=6 C \epsilon$, then $x_1, x_2 \in B(x_0,r_0/2)$. Let $d=
\frac{1}{2}\min\{u(x_1),u(x_2)\}$. We will show that, there exists
$c\geq1$, such that if $\ov{\delta} \leq 1/(48c)$, and $R=c r_0
\leq 1/8$, then the connected components $A_i$ of $B(x_0,R) \cap
B_{1}^d(u)$ which contain $x_i, i=1,2$, are the same. Indeed, let
us suppose that $A_1 \neq A_2$ and let us use Lemmas \ref{monot},
and Corollary \ref{lowbound} with $v=(u-d)^+.$ The density
property of $u$ guarantees that the hypotheses of Lemma
\ref{monot} are satisfied, hence, for some exponent $\beta
>0$, the function $r^{- \beta }\phi(r)$ is non-decreasing. By
Corollary \ref{lowbound} we obtain \[\phi(r_0)
> \tau^2.\]

\noindent Moreover, since $u$ in Lipschitz on $B_{3/4}$ we also
have the bound
\[\phi(R) \leq c',\] with $c'$ absolute constant independent of $R.$
Hence,
\[\tau^2 r_0^{-\beta} < r_0^{-\beta} \phi(r_0) \leq R^{-\beta} \phi(R) \leq
c' R^{-\beta}\]or $R < c'r_0,$ which is a contradiction if we
choose $c=c'$.

\noindent We therefore conclude that $A_1=A_2.$ Since $A_1$ is
open and connected we may find a curve $\Gamma$ inside $A_1$
having $x_1$ and $x_2$ as end point. Denote by $m$ the
non-degeneracy constant of $u$ on $B_{3/4}.$ Then, for each $y \in
\Gamma$ we know that
\[u(y)
>d = \frac{1}{2}\min\{m d(x_1,F(u)),m d(x_2,F(u))\} \geq \frac{1}{2}m\epsilon.\]
Therefore, if $K$ is the Lipschitz constant of $u$ on $B_{3/4}$,
for any $y \in \Gamma$, we have $d(y,F(u))> \frac{1}{2}\frac{m
\epsilon}{K}$. Set $\rho = \frac{1}{2}\frac{m \epsilon}{K},$ so
that if $y \in \Gamma$ and $|x-y|< \rho$ then $u(x)>0.$ Since
\[\Gamma \subset \bigcup_{y \in \Gamma} B(y,\rho)\] we may find a
sequence $y_1,...,y_l$ of points in $\Gamma$ such that $\Gamma
\subset \bigcup_{i=1}^{l} B(y_i,\rho),$ and we may further ask
that no $y$ in $\Gamma$ belong to more than $c(n)$ of the balls
$B(y_i,\rho)$.

\noindent Furthermore, since $\rho = \frac{1}{2}\frac{m
\epsilon}{K}, r_0=6C\epsilon$ and $y_i \in B(x_0,cr_0),$ $l$ must
be bounded by a constant depending only on dimension on $c, C$,
but independent of $x_1,x_2$ or $\epsilon.$ \qed

\

\noindent \textbf{Remark 1.} Notice that Lemma \ref{HC} shows that
the free boundary of any viscosity solution to \eqref{FB1}, which
is Lipschitz continuous, non-degenerate, and satisfies the density
property $(D)$, is locally NTA.

\

\noindent  \textbf{Remark 2.} It follows from the proof of Lemma
\ref{HC} that $M$ and $\ov{\delta}$ depend on the Lipschitz and
non-degeneracy constants of $u$ on $B_{3/4}$.

\

Finally, we recall a fundamental result about NTA domain (see
\cite{JK}).

\begin{thm}\label{dalb}(Dalbherg Boundary Harnack principle) Let $\Omega$ be
an NTA domain, and let $V$ be an open set. For any compact set $G
\subset V$, there exists a constant $C$ such that for all positive
harmonic functions $u$ and $v$ in $\Omega$ that vanish
continuously on $\p \Omega \cap V$, $u(x_0)=v(x_0)$ for some $x_0
\in \Omega \cap G$ implies $C^{-1} u(x) < v(x) < C u(x)$ for all
$x \in G \cap \ov{\Omega}.$
\end{thm}

%\begin{thm}\label{hold}Let $\Omega$ be an NTA domain, and let $V$ be an open set. Let $G$ be a compact
%subset of $V$. There exists a number $\alpha > 0,$ such that for
%all positive harmonic functions $u$ and $v$ in $\Omega$ that
%vanish continuously on $\p \Omega \cap V$, the function
%$u(x)/v(x)$ is H$\ddot{o}$lder continuous of order $\alpha$ on $G
%\cap \ov{\Omega}.$ In particular, for every $y \in G \cap \p
%\Omega,$ $\lim_{x\rightarrow y} (u(x)/v(x))$ exists.
%\end{thm}

\section{Global monotone solutions.}

In this section we prove our global result, that is Theorem
\ref{z}. We start by deriving the following existence result, with
a straightforward limit argument.

\begin{thm}\label{exglobal} Assume that, there exist a strict
smooth subsolution $V_1$ and a strict smooth supersolution $V_2$
to $(\ref{FB1})$ in $\R^n$, such that
\begin{enumerate}
\item $V_1 \leq V_2$ on $\R^n, 0 \in \{V_2>0\} \cap \{V_1=0\}^\circ$;
\item $\p_n V_i >0$ in $\ov{\{V_i>0\}},$ for $i=1,2$.
\end{enumerate}
Then, there exists a global function $u,$ viscosity solution to
$(\ref{FB1})$ in $\R^n$, such that $u$ is monotone increasing in
$\{u>0\}$ in the $x_n$ direction. Moreover $V_1\leq u \leq V_2$,
$u$ is Lipschitz continuous, $(I)$ non-degenerate, and it
satisfies the density property $(D).$
\end{thm}
\noindent\textsc{Proof.} Let $\{R_k\}$ be a sequence of radii,
$R_k \rightarrow +\infty$. Set $u_k := u_{R_k}$, where $u_{R_k}$
is the viscosity solution on $C_{R_k}$ from Theorem \ref{localex}.
Then, by Lemma \ref{Lip}, for any compact subset $D \subset \R^n$,
and sufficiently large $k$, the functions $\{u_k\}$ are uniformly
Lipschitz continuous on $D$. Hence, there exists a function $u:
\R^n \longrightarrow \R^+$, such that (up to a subsequence), $u_k
\rightarrow u$ uniformly on compacts of $\R^n,$ hence $V_1 \leq u
\leq V_2.$ Moreover, $u$ is locally Lipschitz continuous, and
monotone increasing in its positive phase in the $x_n$ direction.
Also, since the $u_k$'s are Lipschitz continuous, $(I)$
non-degenerate, and satisfy the density property $(D)$, with
universal local constants, arguing as in Lemma \ref{blowup}, we
obtain:
\begin{enumerate}\renewcommand{\theenumi}{\alph{enumi}}
\item $\p \{u_k >0 \} \rightarrow \p \{u >0 \}$ in the Hausdorff distance;
\item $\chi_{\{u_k > 0\}} \rightarrow \chi_{\{u > 0\}}$ in $L^1_{loc}$;
\item $\nabla u_k \rightarrow \nabla u$ a.e.
\end{enumerate}
In particular, $u$ is non-degenerate, $(I)$ non-degenerate, and
satisfies the density property $(D)$. Furthermore, $u$ is a
variational solution to (\ref{FB1}), on any compact, and it is
harmonic in its positive phase. A blow-up argument as in Lemma
\ref{FBnotouch1} allows us to conclude that $F(u)$ cannot touch
neither $F(V_1)$ nor $F(V_2)$. Hence, $u$ is a viscosity solution
to (\ref{FB1}) in $\R^n.$\qed

\

In particular, we can also conclude the following:

\begin{cor} $u$ minimizes $J(\cdot, B)$ among all competitors $v \in H^1(B)$, such that $V_1 \leq v
\leq V_2$, and $v=u$ on $\partial B$, for all balls $B \subset
\R^n $.\end{cor}

To conclude the proof of Theorem \ref{z}, we need to prove the
following:

\begin{thm} $F(u)$ is a continuous graph, with a universal modulus of continuity on every compact $K \subset \R^n$.
\end{thm}
\textsc{Proof.} We start by proving that $F(u)$ is a graph.
Assume, by contradiction, that $F(u)$ contains a vertical segment.

\noindent Let $v(x)=u(x-t e_n)$, for some small $t$. Since $u$ is
monotone in the $x_n$ direction, we have  $v \leq u$, and $v < u$
in $\{u>0\}$. Moreover, by the assumption that $F(u)$ contains
vertical segments, we have that $F(u) \cap F(v) $ is non-empty,
for $t$  sufficiently small. Assume, without loss of generality,
that $0 \in F(u) \cap F(v)$. From Lemma \ref{weiss}, we obtain
that $u$ and $v$ blow up around $0$ to functions $U$ and $V$ which
are homogeneous of degree 1. Also, $U$ and $V$ are viscosity
solutions in $\R^n,$ locally Lipschitz continuous, non-degenerate,
and satisfying the density property $(D).$ Hence, according to
Remark 2 following Lemma \ref{HC}, their free boundaries are NTA.
Moreover, $U \geq V$.

\noindent We wish to prove that,
\begin{equation}\label{UlV}U=\lambda V, \ \ \text{on $\{V>0\},$}
\end{equation}for some number $\lambda \geq 1.$ Then, since $U$
and $V$ have the same asymptotic development near regular points,
we get
\begin{equation}\label{UV}U=V \ \ \text{ on
$\{V>0\}$}.\end{equation}

\noindent Towards proving \eqref{UlV}, let us set
$$\lambda = \sup\{t>0 |  U \geq t V \ \ \text{in $\{V>0\}$}\}.$$
Clearly, since $U \geq V,$ then $1 \leq \lambda < +\infty.$ Define
$W$ to be the harmonic function in $B_1 \cap \{V>0\}$ with
boundary data

$$W =
  \begin{cases}
    U & \text{on $\partial B_1 \cap \{V>0\}$}, \\
    0 & \text{on $F(V)$}.
  \end{cases}
$$
By the maximum principle, $$\lambda V \leq W \leq U \ \ \text{on
$B_1 \cap \{V>0\}$}.$$ If at some point $x_0 \in B_{1/2} \cap
\{V>0\}$, we have $(W-\lambda V)(x_0)=0,$ then $W \equiv \lambda
V$ on $B_{1} \cap \{V>0\}.$ Thus, since $W=U$ on $\partial B_1
\cap \{V>0\}$, we have that $U=\lambda V$ at some point $x \in
\partial B_1 \cap \{V>0\}$. Therefore, since $U \geq \lambda V$ in $\{V
>0\},$ we conclude that $ U \equiv \lambda V$ on $\{V>0\}.$

If  at some point $x_0 \in B_{1/2} \cap \{V>0\}$, we have
$W(x_0)/(\lambda V(x_0)) = \delta
>1,$ then the boundary Harnack inequality (Theorem \ref{dalb})
implies that $W-\lambda V \geq c\delta \lambda V$ on $B_{1/2} \cap
\{V>0\}$. Therefore, since $U$ and $\lambda V$ are homogeneous of
degree 1, we conclude that $U \geq (1+c\delta) \lambda V$ on
$\{V>0\},$ which contradicts the definition of $\lambda.$

Since $u > v$ in $\{v>0\}$, a similar comparison argument for $u$
and $v$ allows us to conclude that
\begin{equation}\label{uv} u \geq (1+\epsilon)v \ \ \text{on $\{v>0\}\cap B_{1/2}$},\end{equation}
with $\epsilon$ depending on the ratio of $u$ and $v$ at a fixed
scale. Passing to the blow-up limit, \eqref{uv} produces a
contradiction to \eqref{UV}.

%Hence, if $R_j$ is a sequence of radii such that $R_j \rightarrow
%0 $ as $j \rightarrow +\infty$, then $u_j(x)=u(R_j x)/R_j$ and
%$v_j(x)=v(R_j x)/R_j$ converge uniformly on compacts to the same
%function. In particular, since $v$ is $(I)$ non-degenerate, for
%$\epsilon$ small, there is $0 < r \leq 1/2$, such that
%\begin{equation}\label{max} u-v \leq \epsilon
%\max_{B_r(0)} v \ \ \ \text{on $B_r(0).$}\end{equation}Now, let
%$w$ be the solution to the following Dirichlet
%problem:\begin{equation*}
 % \begin{cases}
  %   \Delta w= 0 & \text{on $B_{1}(0)  \cap \{v>0\}$}, \\
   %  w=u-v & \text{on $(\p B_{1}(0)) \cap \{v>0\}$},\\
   % w=0 & \text{on $B_{1} \cap (\p \{v>0\})$}.\\
  %\end{cases}
%\end{equation*}
%By the maximum principle, $u-v \geq w >0$ in $B_1(0) \cap
%\{v>0\}.$ Hence by the Boundary Harnack inequality, $w \geq C v$
%on $B_{1/2}(0) \cap \{v>0\},$ for some constant $C>0$ depending on
%the ratio of $w$ and $v$ at a fixed scale. On the other hand, by
%(\ref{max}) we get $w \leq c \eps v$ on $B_{r}(0) \cap \{v>0\}$.
%Therefore, we get a contradiction if $\epsilon$ is sufficiently
%small.

Now, let us prove that $F(u)$ has a universal modulus of
continuity on each compact cylinder $K$ (box). We will denote by
$$v_K(\{x_n=\phi_1(x')\},\{x_n=\phi_2(x')\}):= \max_{K'} |\phi_1(x')-\phi_2(x')|,$$ the vertical
distance between two graphs in the $x_n$ direction, with $K'$
being the projection of the box $K$ on the hyperplane $\{x_n=0\}$.

\noindent We want to show that for every $K \subset \R^n,$ and any
$\eps
>0$, there exists a $\delta>0$ such that, if $|\eta| < \delta$,
then any $u$ monotone global minimizer for $J$ among competitors
$v, V_1 \leq v \leq V_2$, satisfies,
\[v_K(\{u =\eta\}, F(u)) < \eps.\]
This, together with the non-degeneracy of $u$, gives the desired
modulus of continuity estimate.

\noindent By contradiction, assume that for some $K \subset \R^n$,
there exist a positive number $\eps$, a sequence $\{\eta_j\},
\eta_j \rightarrow 0,$ as $j \rightarrow +\infty$, and a sequence
of energy minimizing solutions $\{u_j\}$, such that
\begin{equation}\label{j}u_j(x_j + \eps e_n ) < \eta_j,\end{equation} for some $x_j \in F(u_j)
\cap K.$

\noindent The uniform Lipschitz continuity of the $u_j$'s, for $j$
large, implies that (up to a subsequence):\[u_j \rightarrow
\widetilde{u}, \ \ \text{uniformly on compacts},\] and \[x_j
\rightarrow \ov{x} \in K,\] with $\widetilde{u}$ a Lipschitz
continuous minimizing solution, monotone increasing in the $x_n$
direction, and satisfying the $(I)$ non-degeneracy, and the
density property $(D)$. Moreover $\widetilde{u}(\ov{x})
=\widetilde{u}(\ov{x} + \eps e_n)=0.$ We aim to prove that $\ov{x}
\in F(\widetilde{u})$; then by (\ref{j}), we obtain that
$F(\widetilde{u})$ contains the vertical segment from $\ov{x}$ to
$\ov{x}+\eps e_n$, which is a contradiction to what we showed
above. Indeed, assume, that $\ov{x}$ does not belong to
$F(\widetilde{u}).$ Then, there exists $r> 0$ such that,
\begin{equation}\label{Br}B_r(\ov{x}) \subset \{\widetilde{u}=0\}^\circ.\end{equation} If $j$ is
large enough, then $x_j \in B_{r/4}(\ov{x}) \cap F(u_j)$, and by
non-degeneracy $$\fint_{B_{r/2}(\ov{x})} u_j \geq K r.$$ Hence we
get a contradiction to \eqref{Br}, as we pass to the limit for $j
\rightarrow \infty$. \qed

\subsection*{Acknowledgements.} The results of this paper are part
of my doctoral dissertation \cite{D} at MIT. I would like to thank
my advisor D. Jerison for introducing me to this subject, and for
his guidance during my stay at MIT.

% ----------------------------------------------------------------
%\bibliographystyle{amsplain}
%\bibliography{}

\begin{thebibliography}{9999}\bibitem[AC]{AC} H.W. Alt, L.A. Caffarelli,
\emph{Existence and regularity for a minimum problem with free
boundary}, J. Reine Angew. Math {\bf 325} (1981),105--144.
\bibitem[ACF]{ACF} H.W. Alt, L.A. Caffarelli, A. Friedman,
\emph{Variational problems with two phases and their free boundary,}
Trans. Amer. Math. Soc.  \textbf{282}  (1984),  no. 2, 431--461.
\bibitem[ACS]{ACS} N.E. Aguilera, L.A. Caffarelli, J.
Spruck, \emph{An optimization problem in heat conduction},  Ann.
Scuola Norm. Sup. Pisa Cl. Sci. (4) \textbf{14} (1987), no. 3,
355--387.
\bibitem[AH]{AH} D.R. Adams, L.I. Hedberg, \emph{Function spaces and potential theory},
Grundlehren der Mathematischen Wissenschaften [Fundamental
Principles of Mathematical Sciences], \textbf{314}
Springer-Verlag, Berlin, (1996) \bibitem[BDG]{BDG} E. Bombieri, E.
De Giorgi, E. Giusti, \emph{Minimal cones and the Bernstein
problem}, Inv. Math. \textbf{7} (1969), 243--268.\bibitem[C1]{C1}
L.A. Caffarelli, \emph{A Harnack inequality approach to the
regularity of free boundaries. Part I: Lipschitz free boundaries
are $C^{1,\alpha}$}, Rev. Mat. Iberoamericana \textbf{3} (1987)
no. 2, 139--162.
\bibitem[C2]{C2} L.A. Caffarelli, \emph{A Harnack
inequality approach to the regularity of free boundaries. Part II:
Flat free boundaries are Lipschitz}, Comm. Pure Appl. Math.
\textbf{42} (1989), no.1, 55--78.\bibitem[C3]{C3} L.A. Caffarelli, \emph{A
Harnack inequality approach to the regularity of free boundaries.
Part III: Existence theory, compactness, and dependence on $X$},
Ann. Scuola Norm. Sup. Pisa Cl. Sci. (4) \textbf{15}
(1988), no. 4, 583--602 (1989).\bibitem[CJK]{CJK} L.A. Caffarelli, D. Jerison, C.E. Kenig,
\emph{Global energy minimizers for free boundary problems and full
regularity in three dimension}, Noncompact problems at the
intersection of geometry, analysis, and topology, 83--97, Contemp.
Math., 350, Amer. Math. Soc., Providence, RI, 2004.
\bibitem[D]{D} D. De Silva, \emph{Existence and
regularity of monotone solutions to free boundary problems}, Ph.D
Thesis, 2005
\bibitem[DJ]{DJ} D. De Silva, D. Jerison, \emph{A singular energy
minimizing free boundary}, Preprint (2005), submitted for
publication.
\bibitem[DJ2]{DJ2} D. De Silva , D. Jerison , \emph{Gradient bound for energy minimizing free
boundary graphs}, in preparation.
\bibitem[F]{F} A. Friedman,
\emph{Variational principles and free-boundary problems}, R.
Krieger Publishing Co., Inc.,Malabar, FL, 1988.
\bibitem[G]{G} E. Giusti,
\emph{Minimal surfaces and functions of bounded variation},
Monographs in Mathematics, 80. Birkh$\ddot{a}$user Verlag, Basel, 1984.
\bibitem[GT]{GT} D. Gilbarg, N.
Trudinger, \emph{Elliptic partial differential equations of second
order,} Grundlehren der Mathematischen Wissenschaften, Vol. 224. Springer-Verlag,
Berlin-New York, 1977.
\bibitem[JK]{JK} D. Jerison, C.E.
Kenig, \emph{Boundary behavior of harmonic functions in
Non-tangentially Accessible Domains}, Adv. in Math. \textbf{46}
(1982), no. 1, 80--147.
\bibitem[K]{K} B. Kawohl, \emph{Rearrangements and Convexity
of Level Sets in PDE, } Lecture Notes in Mathematics, \textbf{1150}.
Springer-Verlag, Berlin, 1985.
\bibitem[KN]{KN} D. Kinderlehrer, L. Nirenberg , {\it Analyticity at the boundary of solutions of
nonlinear second-order parabolic equations,} Comm. Pure Appl.
Math.  31  (1978), no. 3, 283--338.
\bibitem[W1]{W1} G.S. Weiss, \emph{Partial regularity for weak
solutions of an elliptic free boundary problem,} Commun. in
Partial Diff. Equat. \textbf{23} (1998), no. 3-4, 439-455. \bibitem[W2]{W2} G.S.
Weiss, \emph{Partial regularity for a minimum problem with free
boundary,} J. Geom. Anal.  9  (1999),  no. 2, 317--326.
\end{thebibliography}
%\include{biblio}
\end{document}